\documentclass{amsart}

\usepackage{mathrsfs}
\usepackage{amsfonts}
\usepackage{amssymb,amsmath}
\usepackage{amsthm}

\usepackage{graphics}
\usepackage{graphicx}
\usepackage{epsfig}

\newtheorem{theorem}{Theorem}[section]
\newtheorem{lemma}[theorem]{Lemma}
\newtheorem{corollary}[theorem]{Corollary}

\theoremstyle{definition}

\newtheorem{proposition}[theorem]{Proposition}

\theoremstyle{remark}
\newtheorem{remark}[theorem]{Remark}
\theoremstyle{assumption}

\theoremstyle{assumptions}
\newtheorem{assumptions}[theorem]{Assumptions}

\numberwithin{equation}{section}

\usepackage{changes}

\usepackage{setspace}
\newenvironment{citeRef}
{\vskip 0.1cm \begin{singlespace*}\textbf{References:}\begin{small}}
		{\end{small}\end{singlespace*}}

\DeclareMathOperator*{\argmin}{\text{argmin}}
\usepackage[]{algorithm}
\usepackage{algorithmic} 
\usepackage{makecell}
\usepackage{environ}
\NewEnviron{myeq}{
	\begin{align}
	\begin{split}
	\BODY 
	\end{split}
	\end{align}
}
\newcommand{\mye}[1]{
	\begin{myeq}
		#1
	\end{myeq}
}
\NewEnviron{myeqn}{
	\begin{align*}
	\BODY
	\end{align*}
}
\newcommand{\myen}[1]{
	\begin{myeqn}
		#1
	\end{myeqn}
}
\NewEnviron{myeql}{
	\begin{align}\left\{\begin{aligned}
	\BODY
	\end{aligned}\right.\end{align}
}

\NewEnviron{myeqln}{
	\begin{align*}\left\{\begin{aligned}
	\BODY
	\end{aligned}\right.\end{align*}
}

\makeatletter
\newcommand{\xRightarrow}[2][]{\ext@arrow 0359\Rightarrowfill@{#1}{#2}}
\makeatother

\begin{document}

\title[Variational inference for functions]
 {Variational Bayes' method for functions with applications to some inverse problems}

\author[J. Jia]{Junxiong Jia}
\address{School of Mathematics and Statistics,
Xi'an Jiaotong University,
Xi'an,
710049, China}
\email{jjx323@xjtu.edu.cn}

\author[Q. Zhao]{Qian Zhao}
\address{School of Mathematics and Statistics,
	Xi'an Jiaotong University,
	Xi'an,
	710049, China}
\email{timmy.zhaoqian@xjtu.edu.cn}

\author[Z. Xu]{Zongben Xu}
\address{School of Mathematics and Statistics,
	Xi'an Jiaotong University,
	Xi'an,
	710049, China}
\email{zbxu@mail.xjtu.edu.cn}

\author[D. Meng]{Deyu Meng*}
\address{School of Mathematics and Statistics,
	Xi'an Jiaotong University,
	Xi'an,
	710049, China}
\email{dymeng@mail.xjtu.edu.cn}

\author[Y. Leung]{Yee Leung}
\address{The Chinese University of Hong Kong, Shatin, Hong Kong, SAR}
\email{yeeleung@cuhk.edu.hk}






	\maketitle
	
	\begin{abstract}
		Bayesian approach, a useful tool for quantifying uncertainties, has been extensively employed to solve the inverse problems of partial differential equations (PDEs). One of the main difficulties in employing the Bayesian approach to such problems is
		how to extract information from the posterior probability measure.
		Compared with conventional sampling-type methods,
		variational Bayes' method (VBM) has been intensively examined in the field of machine learning attributed to its ability in extracting
		approximately the posterior information with lower computational cost.
		In this paper, we generalize the conventional finite-dimensional VBM to the infinite-dimensional space rigorously solve the
		inverse problems of PDEs. We further establish a general infinite-dimensional mean-field approximate theory and apply it to
		the linear inverse problems under the Gaussian and Laplace noise assumptions at the abstract level.
		The results of some numerical experiments substantiate the effectiveness of the proposed approach.
	\end{abstract}
	
	\maketitle
	
	
	\section{Introduction}\label{introduction}
	Motivated by the significant applications in medical imaging, seismic explorations and many other domains,
	the field of inverse problems has undergone an enormous development over the past few decades.
	In handling an inverse problem, we usually meet ill-posed issue in the sense that
	the solution lacks stability or even uniqueness \cite{Ito2015Book,Tarantola2005Book}.
	The regularization approach, including Tikhonov and Total-Variation regularization,
	is one of the most popular approaches to alleviate this ill-posed issue of inverse problems.
	Regularization method allows for statistical assumptions of data, however there are no statistical
	descriptions of the model parameters. For a complete review, we refer to Sections 2 and 3 in \cite{Arridge2019ActaNumerica}. 
	
	The Bayesian inverse approach provides a flexible framework that solves inverse problems by transforming them into
	statistical inference problems, thereby making it feasible to analyze the uncertainty of the solutions to the inverse problems.
	Inverse problems are usually accompanied by a forward operator originating from some partial differential equations (PDEs),
	thereby introducing difficulties to the direct use of the finite-dimensional Bayes' formula.
	The following two strategies can be employed to solve this problem:
	\begin{enumerate}
		\item Discretize-then-Bayesianize: The PDEs are initially discretized to approximate the original problem in
		some finite-dimensional space, and the reduced approximate problem is then solved by using the Bayes' method.
		\item Bayesianize-then-discretize: The Bayes' formula and algorithms are initially constructed on infinite-dimensional space,
		and after the infinite-dimensional algorithm is built, some finite-dimensional approximation is carried out.
	\end{enumerate}
	The first strategy makes available all the Bayesian inference methods developed in the statistical literature \cite{Kaipio2005Book}.
	However, given that the original problems are defined on infinite-dimensional space, several problems, such as non-convergence and
	dimensional dependence, tend to emerge when using this strategy \cite{Cotter2013SS, Lassas2004IP}. By employing the second strategy,
	the discretization-invariant property naturally holds given that the Bayes' formula and algorithms are properly defined on some
	separable Banach space \cite{Dashti2014,Stuart2010AN}. In the following sections, we confine ourselves to the second strategy, that is,
	postponing the discretization to the final step.
	
	One of the essential issues for employing the Bayes' inverse method is how to extract
	information from the posterior probability measure. Previous studies have adopted two
	major approaches to address such issue, namely, the point estimate method and the sampling method.
	For the point estimate method, the maximum a posteriori (MAP) estimate, which is intuitively
	equivalent to solving an optimization problem, is often utilized.
	The intuitive equivalence relation has been rigorously analyzed recently \cite{Agapiou2018IP,Burger2014IP,Dashti2013IP,Dunlop2016IP,Helin2015IP}.
	In some situations \cite{Jia2019IP,Tarantola2005Book}, MAP estimates are more desirable and computationally feasible than
	the entire posterior distribution. However, point estimates cannot provide uncertainty quantification and are usually
	recognized as incomplete Bayes' method.
	
	To extract all information encoded in the posterior distribution, sampling methods, such as the Markov chain Monte Carlo (MCMC),
	are often employed. In 2013, Cotter et al. \cite{Cotter2013SS} proposed using the MCMC method for functions to
	ensure that the convergence speed of the algorithm is robust under mesh refinement.
	Multiple dimension-independent MCMC-type algorithms have also been proposed \cite{Cui2016JCP,Feng2018SIAM}.
	Although MCMC is highly-efficient as a sampling method, its computational cost
	is unacceptable for many applications, including the full waveform inversion \cite{Fichtner2011Book}.
	
	In this paper, we aim to propose a variational method that can perform uncertainty analysis at a
	computational cost which is comparable to that for computing the MAP estimates. For finite-dimensional problems,
	such types of methods, named as variational Bayes' methods (VBM), have been broadly investigated
	in the field of machine learning \cite{Bishop2006PRML,Matthews2016PHD,Zhang2018IEEE,Zhao2015IEEE}.
	In addressing the inverse problems, Jin et al. \cite{Jin2010JCP,Jin2012JCP} employed VBM to investigate
	a hierarchical formulation of the finite-dimensional inverse problems when the noise is distributed according
	to Gaussian or centered-t distribution.
	Guhua et al. \cite{Guha2015JCP} generalized this method further to the case when the noise is distributed according
	to skewed-t error distribution. Finite-dimensional VBM has been recently applied
	to study the porous media flows in heterogeneous stochastic media \cite{Yang2017JCP}.
	
	All the aforementioned investigations are conducted based on finite-dimensional VBM.
	Therefore, only the first strategy as aforementioned can be employed to solve the inverse problems.
	To the best of our knowledge, only two relevant works have investigated VBM under the infinite-dimensional setting.
	Specifically, when the approximate probability measures are restricted to be Gaussian,
	Pinski et al. \cite{Pinski2015SIAMSC,Pinski2015SIAMMA} employed a calculus-of-variations viewpoint to study
	the properties of Gaussian approximate sequences with Kullback-Leibler (KL) divergence as the a fitness measure.
	Relying on the Robbins-Monro algorithm, they developed a novel algorithm for obtaining the approximate Gaussian measure.
	Until now, no study has been conducted beyond such Gaussian approximate measure assumption.
	However, various approximate probability measures have been frequently used for training deep neural networks and
	solving finite-dimensional inverse problems \cite{Jin2010JCP,Jin2012JCP}.
	In this case, for applications in inverse problems concerned with PDEs,
	a VBM with approximate measures other than Gaussian should be necessarily constructed on infinite-dimensional space.
	
	In the following, we focus on the classical mean-field approximation that is widely employed for the finite-dimensional case.
	This approximation originally stems from the theory of statistical mechanics for treating many-body systems.
	Inspired by finite-dimensional theory, we construct a general infinite-dimensional mean-field approximate based VBM,
	which allows the use of general approximate probability measures beyond Gaussian.
	Examples are also given to illustrate the flexibility of our proposed approach.
	The contributions of our work can be summarized as follows:
	\begin{itemize}
		\item By introducing a reference probability measure and using the calculus of variations,
		we establish a general mean-field approximate based VBM on Hilbert spaces that provides a
		flexible framework for introducing techniques developed on finite-dimensional space to infinite-dimensional space.
		\item We apply the proposed theory to a general linear inverse problem
		(the forward map is assumed to be a bounded linear operator)
		with Gaussian and Laplace noise assumptions.
		Precise assumptions can be found in Subsection \ref{SectionLinearIPG}.
		Through detailed calculations, we construct iterative algorithms for functions.
		To the best of our knowledge, VBM with Laplace noise assumption has not been previously employed for solving inverse problems,
		even those that are restricted to finite-dimensional space.
		\item We solve the inverse source problems of Helmholtz equations with multi-frequency data by using the proposed VBM
		with Gaussian and Laplace noise assumptions.	
		The algorithms not only provide a point estimate but also give the standard deviations of the numerical solutions.
	\end{itemize}
	
	The outline of this paper is as follows.
	In Section 2, we construct the general infinite-dimensional VBM based on the mean-field approximate assumption.
	In Section 3, under the hierarchical formulation, we apply the proposed theory to an abstract linear inverse problem
	with Gaussian and Laplace noise assumptions.
	In Section 4, we present concrete numerical examples to illustrate the effectiveness of our proposed approach.
	In Section 5, we summarize our findings and propose some directions for further research.
	Due to the limited space, we did not provide all proofs in the main text.
	\emph{All of the proofs are given in the supplemental materials}.
	
	\section{General theory on infinite-dimensional space}\label{SectionGeneralTheory}
	
	In Subsection \ref{ExistenceSection}, we provide the necessary background of our theory and
	prove some basic results concerning with the existence of minimizers for finite product probability measures.
	In Subsection \ref{MeanFieldSection}, we present our infinite-dimensional variational Bayes' approach.
	
	\subsection{Existence theory}\label{ExistenceSection}
	In this subsection, we first recall some general facts about the Kullback-Leibler (KL) approximation from the viewpoint of
	calculus of variations, and then provide some new theorems for product of probability measures that form
	the basis of our investigation. Let $\mathcal{H}$ be a Polish space endowed with its Borel sigma algebra
	$\mathcal{B}(\mathcal{H})$, and	let $\mathcal{M}(\mathcal{H})$ be the set of Borel probability measures
	on $\mathcal{H}$.
	
	For inverse problems, we usually need to find a probability measure $\mu$ on $\mathcal{H}$,
	which is called the posterior probability measure, specified by its density with respect to a prior
	probability measure $\mu_{0}$ \cite{Stuart2010AN}. Let the Bayesian formula on the Hilbert space be defined by
	\mye{\label{BayesianFormula1}
		\frac{d\mu}{d\mu_{0}}(x) = \frac{1}{Z_{\mu}}\exp\big( -\Phi(x) \big),
	}
	where $\Phi(x): \mathcal{H}\rightarrow \mathbb{R}$ is a continuous function, and
	$\exp\big( -\Phi(x) \big)$ is integrable with respect to $\mu_{0}$. The constant $Z_{\mu}$ is chosen
	to ensure that $\mu$ is indeed a probability measure.
	
	Let $\mathcal{A}\subset\mathcal{M}(\mathcal{H})$ be a set of ``simpler'' measures that can be efficiently calculated.
	Our aim is to find the closest element $\nu$ to $\mu$ with respect to the KL divergence from
	subset $\mathcal{A}$. For any $\nu\in\mathcal{M}(\mathcal{H})$ that is absolutely continuous with respect to $\mu$,
	the KL divergence is defined as
	\mye{
		D_{\text{KL}}(\nu||\mu) =
		\int_{\mathcal{H}}\log\bigg( \frac{d\nu}{d\mu}(x) \bigg)\frac{d\nu}{d\mu}(x)\mu(dx)
		= \mathbb{E}^{\mu}\bigg[ \log\bigg( \frac{d\nu}{d\mu}(x) \bigg)\frac{d\nu}{d\mu}(x) \bigg],
	}
	where the convention $0\log 0 = 0$ has been used. If $\nu$ is not absolutely continuous with respect to $\mu$,
	then the KL divergence is defined as $+\infty$.
	With this definition, this paper examines the following minimization problem:
	\mye{\label{geneMiniPro1}
		\argmin_{\nu\in\mathcal{A}}D_{\text{KL}}(\nu||\mu).
	}
	
	There are some studies of the above general minimization problem (\ref{geneMiniPro1}) taken from the perspective of the calculus of variations.
	We follow this line of investigations in this section, and for the convenience of the readers, we present
	the following proposition, which has been proven in \cite{Pinski2015SIAMMA}.
	
	\begin{proposition}\label{por1}
		Let $\mathcal{A}$ be closed with respect to weak convergence. Then, given $\mu\in\mathcal{M}(\mathcal{H})$,
		assume that there exists $\nu\in\mathcal{A}$ such that $D_{KL}(\nu||\mu) < \infty$.
		It follows that there exists a minimizer $\nu\in\mathcal{A}$ solving
		\begin{align*}
		\argmin_{\nu\in\mathcal{A}}D_{\text{KL}}(\nu||\mu).
		\end{align*}
	\end{proposition}
	
	As stated in the Introduction, we aim to construct a mean-field approximation that
	usually takes the following factorized form for the finite-dimensional case
	\mye{\label{fineDmeanField1}
		q(x_{1},\cdots,x_{M}) = \prod_{j=1}^{M}q(x_{j}),
	}
	where $q(x_1,\cdots,x_{M})$ is the full probability density function, $q(x_{j})$ is the probability density function for $x_{j}$,
	and $x_{j}\in\mathbb{R}^{N_{j}}$($N_{j}\in\mathbb{N}$) for $j=1,2,\cdots,M$.
	That is, we assume that $x_1,\cdots,x_M$ are independent random variables.
	By carefully choosing the random variables $\{x_{j}\}_{j=1}^{M}$, this independence assumption will lead to
	computationally efficient solutions when conjugate prior probabilities are employed.
	Additional details can be found in Chapter 9 of \cite{Bishop2006PRML} and in some recently published papers
	\cite{Jin2012JCP,Zhang2018IEEE,Zhao2015IEEE}.
	
	Inspired by formula (\ref{fineDmeanField1}), for a fixed positive constant $M$, 
	we specify the Hilbert space $\mathcal{H}$ and subset $\mathcal{A}$
	as
	\mye{\label{prodHA}
		\mathcal{H} = \prod_{j=1}^{M}\mathcal{H}_{j}, \qquad
		\mathcal{A} = \prod_{j=1}^{M}\mathcal{A}_{j},
	}
	where $\mathcal{H}_{j} (j=1,\cdots,M)$ are a series of separable Hilbert space and
	$\mathcal{A}_{j}\subset\mathcal{M}(\mathcal{H}_{j})$.
	Let $\nu:=\prod_{i=1}^{M}\nu^{i}$ be a probability measure such that $\nu(dx) = \prod_{i=1}^{M}\nu^{i}(dx)$.
	With these assumptions, the minimization problem in (\ref{geneMiniPro1}) can be rewritten as
	\mye{\label{mipro1}
		\argmin_{\nu^{i}\in\mathcal{A}_{i}}D_{\text{KL}}\bigg(\prod_{i=1}^{M}\nu^{i}\big|\big|\mu\bigg)
	}
	for suitable sets $\mathcal{A}_{i}$ with $i=1,2,\cdots,M$.
	The general result shown in Proposition \ref{por1} indicates that the closedness of the subset $\mathcal{A}$
	under weak convergence is crucial for the existence of the approximate measure $\nu$. Therefore, we
	present the following lemma that illustrates the closedness of $\mathcal{A}$ as defined in (\ref{prodHA}).
	
	\begin{lemma}\label{sepLe1}
		For $i=1,2,\cdots, M$, let $\mathcal{A}_{i} \subset \mathcal{M}(\mathcal{H}_{i})$ be a series of sets
		closed under weak convergence of probability measures. Define
		\mye{
			\mathcal{C}:= \bigg\{ \nu:=\prod_{j=1}^{M}\nu^{j} \,\, \bigg| \,\,
			\nu^{j}\in\mathcal{A}_{j} \quad \text{for}\,\,\, j=1,2,\cdots,M \bigg\}.
		}
		Then, the set $\mathcal{C}$ is closed under the weak convergence of probability measures.
	\end{lemma}
	
	From Lemma \ref{sepLe1} and Proposition \ref{por1}, we can prove the following existence result.
	\begin{theorem}\label{genetheo1}
		For $i=1,2,\cdots,M$, let $\mathcal{A}_{i}$ be closed with respect to weak convergence.
		Given $\mu\in\mathcal{M}(\prod_{i=1}^{M}\mathcal{H}_{i})$, we assume that there exists $\nu^{i}\in\mathcal{A}_{i}$
		for $i=1,\cdots,M$
		such that $D_{\text{KL}}(\prod_{i=1}^{M}\nu^{i} || \mu)\! < \!\infty$. \hspace{-0.8em}
		Then, there exists a minimizer $\prod_{i=1}^{M}\nu^{i}$ that solves problem (\ref{mipro1}).
	\end{theorem}
	
	\begin{remark}\label{noUniquenessRemark}
		In Theorem \ref{genetheo1}, we only illustrate the existence of the approximate measure $\nu$ without uniqueness.
		When the approximate measures are assumed to be Gaussian, uniqueness has been obtained with the
		$\lambda$-convex requirement of the potential $\Phi$ appearing in
		the Bayes' formula (\ref{BayesianFormula1}) \cite{Pinski2015SIAMMA}.
		We cannot expect uniqueness generally even for most of the practical problems defined on the finite-dimensional space.
		Therefore, we will not pursue the uniqueness results here.
	\end{remark}
	
	The result shown in Theorem \ref{genetheo1} does not tell us much about the manner in which
	minimizing sequences approach the limit. After further deductions, we can precisely characterize
	the convergence.
	
	\begin{theorem}\label{genetheo2}
		Let $\big\{\nu_{n}=\prod_{j=1}^{M}\nu_{n}^{j}\big\}_{n=1}^{\infty}$ be a sequence in
		$\prod_{j=1}^{M}\mathcal{M}(\mathcal{H}_{j})$, and let
		$\nu_{*}=\prod_{j=1}^{M}\nu_{*}^{j} \in \prod_{j=1}^{M}\mathcal{M}(\mathcal{H}_{j})$
		and $\mu \in \mathcal{M}(\prod_{j=1}^{M}\mathcal{H}_{j})$
		be probability measures such that for any $n\geq 1$, we have $D_{\text{KL}}(\nu_{n}||\mu) < \infty$ and
		$D_{\text{KL}}(\nu_{*}||\mu) < \infty$. Suppose that $\nu_{n}$ converges weakly to $\nu_{*}$ and
		$\nu_{n}^{j} \ll \nu^{j}_{*}$ for $j=1,2,\cdots, M$ and that
		\mye{
			D_{\text{KL}}(\nu_{n}||\mu) \rightarrow D_{\text{KL}}(\nu_{*}||\mu).
		}
		Then, $\nu_{n}^{j}$ converges to $\nu_{*}^{j}$ in the total variation norm for $j=1,2,\cdots,M$.
	\end{theorem}
	
	Combining Theorems \ref{genetheo1} and \ref{genetheo2}, we immediately obtain the following result.
	\begin{corollary}\label{genecorollary1}
		For $j=1,2,\cdots, M$, let $\mathcal{A}_{j}\subset\mathcal{M}(\mathcal{H}_{j})$ be closed with respect to
		weak convergence. Given $\mu\in\mathcal{M}(\prod_{j=1}^{M}\mathcal{H}_{j})$, there exists
		$\nu = \prod_{j=1}^{M}\nu^{j}\in\prod_{j=1}^{M}\mathcal{A}_{j}$ with $D_{\text{KL}}(\nu||\mu) <\infty$.
		Let $\nu_{n}=\prod_{j=1}^{M}\nu_{n}^{j}\in\prod_{j=1}^{M}\mathcal{A}_{j}$ satisfy
		\mye{\label{gr1}
			D_{\text{KL}}(\nu_{n}||\mu) \rightarrow \inf_{\nu\in\prod_{j=1}^{M}\mathcal{A}_{j}}D_{\text{KL}}(\nu||\mu).
		}
		Then, after passing to a subsequence, we have
		\begin{itemize}
			\item $\nu_{n}$ converges weakly to $\nu_{*}=\prod_{j=1}^{M}\nu_{*}^{j}\in\prod_{j=1}^{M}\mathcal{M}(\mathcal{H}_{j})$
			that realizes the infimum in (\ref{gr1});
			\item each $\nu_{n}^{j}$ converges weakly to $\nu_{*}^{j}$ for $j=1,2,\cdots,M$.
		\end{itemize}
		In addition, for $j=1,2,\cdots,M$, if $\nu_{n}^{j} \ll \nu_{*}^{j}$ for all $n$,
		each $\nu_{n}^{j}$ converges to $\nu_{*}^{j}$ in the total-variation norm.
	\end{corollary}
	
	
	\subsection{Mean-field approximation for functions}\label{MeanFieldSection}
	
	For finite-dimensional cases, the mean-field approximation has been widely employed in the field of
	machine learning. On the basis of the results presented in Subsection \ref{ExistenceSection}, we construct
	a mean-field approximation approach on infinite-dimensional space, which will be useful for solving the inverse problems of PDEs.
	
	In the previous work, e.g., Examples 3.8 and 3.9 in \cite{Pinski2015SIAMMA} and the general setting described in \cite{Pinski2015SIAMSC},
	their idea is replacing the classical density functions by the density functions with respect to the prior measure.
	In \cite{Pinski2015SIAMSC,Pinski2015SIAMMA}, prior measures are taken to be Gaussian measures, which plays the role of Lebesgue measure
	in the finite-dimensional setting. Inspired by these studies, we may assume that the approximate probability measure $\nu$
	introduced in (\ref{geneMiniPro1}) is equivalent to $\mu_0$ defined by
	\mye{
		\frac{d\nu}{d\mu_{0}}(x) = \frac{1}{Z_{\nu}}\exp\big( -\Phi_{\nu}(x) \big).
	}
	Compared with the finite-dimensional case, a natural way for introducing an independence assumption is to
	assume that the potential $\Phi_{\nu}(x)$ can be decomposed as
	\mye{
		\exp\left(-\Phi_{\nu}(x)\right) = \prod_{j=1}^{M}\exp\left(-\Phi_{\nu}^{j}(x_{j})\right),
	}
	where $x = (x_{1},\cdots,x_{M})$. However, this intuitive idea prevents us from incorporating those parameters contained
	in the prior probability measure into the hierarchical Bayes' model that is used in
	finite-dimensional cases \cite{Jin2012JCP,Zhao2015IEEE}.
	Given these considerations, we propose the following assumption that introduces a reference probability measure.
	
	\begin{assumptions}\label{assumption1}
		Let us introduce a reference probability measure
		\mye{\mu_{r}(dx) = \prod_{j=1}^{M}\mu_{r}^{j}(dx_{j}),}
		which is equivalent
		to the prior probability measure with the following relation being true:
		\mye{
			\frac{d\mu_{0}}{d\mu_{r}}(x) = \frac{1}{Z_{0}}\exp(-\Phi^{0}(x)).
		}
		For each $j=1,2,\cdots,M$, there is a predefined continuous function $a_j(\epsilon, x_j)$\footnote{These functions naturally appear when considering concrete examples, which will be specified in Remark \ref{remarkEx1}.
			In the last part of the supplementary materials, we provide a detailed illustration of the Gaussian noise example, which may provide more intuitions.}
		where $\epsilon$ is a positive number and $x_j\in\mathcal{H}_j$.
		Concerning these functions, we assume that $\mathbb{E}^{\mu_r^j}[a_j(\epsilon,\cdot)] < \infty$
		where $\epsilon\in[0, \epsilon_0^j)$ with $\epsilon_0^j$ is a small positive number
		($j=1,\cdots,M$).
		We also assume that the approximate probability measure $\nu$ is equivalent to the reference measure $\mu_{r}$
		and that the Radon-Nikodym derivative of $\nu$ with respect to $\mu_{r}$ takes the following form
		\mye{\label{defP2}
			\frac{d\nu}{d\mu_{r}}(x) = \frac{1}{Z_{r}}\exp\bigg( -\sum_{j=1}^{M}\Phi_{j}^{r}(x_{j}) \bigg).
		}
	\end{assumptions}
	
	Following Assumptions \ref{assumption1}, we know that the approximate measure can be decomposed as
	$\nu(dx) = \prod_{j=1}^{M}\nu^{j}(dx_j)$ with
	\mye{\label{aaa1}
		\frac{d\nu^{j}}{d\mu_{r}^{j}} = \frac{1}{Z_{r}^{j}}\exp\big( -\Phi_{j}^{r}(x_{j}) \big).
	}
	Here, $Z_{r}^{j} = \mathbb{E}^{\mu_{r}^{j}}\big(\exp\big( -\Phi_{j}^{r}(x_{j}) \big)\big)$ ensures that $\nu^{j}$
	is indeed a probability measure.
	
	\begin{remark}
		The reference measure introduced above can be easily specified for concrete examples. Fix a component $j$,
		if $x_j$ belongs to some finite-dimensional Hilbert space, we assume that the prior measure of $x_j$ has a density function $p(\cdot)$.
		Then we can choose the reference measure of $x_j$ just equal to the prior measure.
		Formula (\ref{aaa1}) for this component reduces to the classical finite-dimensional case.
		If $x_j$ belongs to some Hilbert space with the prior measure contains some hyper-parameters,
		there may be no universal strategies for choosing the reference measure.
		Here, we provide a simple example to give some intuitions.
		Assume $x_j\sim\mathcal{N}(0,\mathcal{C}_{\tau})$ with $\mathcal{C}_{\tau}:=(\tau^2 I-\Delta)^{-\alpha}$
		($\alpha$ is a fixed positive number) \cite{Dunlop2017SC}, we can choose the reference measure to be a Gaussian measure $\mathcal{N}(0,\mathcal{C})$
		with $\mathcal{C}:=(I-\Delta)^{-\alpha}$, which is equivalent to $\mathcal{N}(0,\mathcal{C}_{\tau})$ under some appropriate conditions
		(rigorous results are given in Theorem 1 in \cite{Dunlop2017SC}).
	\end{remark}
	
	For convenience, let us introduce some notations. For $j=1, 2,\cdots, M$, let $\mathcal{Z}_j$ be defined
	as a Hilbert space that is embedded in $\mathcal{H}_j$.
	Denote $C_N$ be a positive constant related $N$.
	Then, for $j=1,2,\cdots,M$, we introduce
	\begin{align*}
	& \qquad\qquad\,\,
	\text{R}^1_j = \bigg\{
	\Phi_j^r \, \Big| \, \sup_{1/N \leq \|x_j\|_{\mathcal{Z}_j} \leq N} \Phi_j^r(x_j)\leq C_N < \infty
	\text{ for all }N>0
	\bigg\},   \\
	& \text{R}^2_j =
	\bigg\{\Phi_j^r \, \Big| \, \int_{\mathcal{H}_j}\exp\left( -\Phi_{j}^{r}(x_j) \right)
	\max(1, a_j(\epsilon,x_j)) \mu_r^{j}(dx_j) \leq C < \infty, \text{ for }\epsilon\in[0,\epsilon_0^j) \bigg\},  
	\end{align*}
	where $C$ is an arbitrary large positive constant,
	$\epsilon_0^j$ and $a_j(\cdot,\cdot)$ are defined as in Assumptions \ref{assumption1}.
	With these preparations, we can define $\mathcal{A}_{j}$ ($j=1,2,\cdots,M$) as follows:
	\begin{align}\label{defAj}
	\mathcal{A}_{j} = \bigg\{ \nu^{j}\in\mathcal{M}(\mathcal{H}_{j}) \, \bigg| \,
	\begin{array}{l}
	\nu^{j} \text{ is equivalent to } \mu_{r}^{j} \text{ with }
	(\ref{aaa1}) \text{ holds true, }  \\
	\text{and }
	\Phi_{j}^{r}\in \text{R}^1_j \cap \text{R}^2_j
	\end{array}
	\bigg\}.
	\end{align}
	Before using Theorem \ref{genetheo1}, we need to illustrate the closedness of $\mathcal{A}_{j}$
	($j=1,2,\cdots,M$) under the weak convergence topology.
	Actually, we can prove the desired results shown blow.
	\begin{theorem}\label{generalAiClosed}
		For $j=1,2,\cdots,M$,
		we denote $T_{N}^{j} = \{ x_j \, |\, 1/N\leq \|x_j\|_{\mathcal{Z}_{j}} \leq N \}$, with $N$ being an arbitrary positive
		constant. For each reference measure $\mu_{r}^{j}$, we assume that $\sup_{N} \mu_{r}^{j}(T_{N}^{j}) = 1$.
		Then, each $\mathcal{A}_{j}$ is closed with respect to weak convergence
		and problem (\ref{mipro1}) possesses a solution
		$\prod_{j=1}^{M}\nu^{j}$ with $\nu^j \in \mathcal{A}_j$ for $j=1,2,\cdots, M$.
	\end{theorem}
	
	In the following theorem, we provide a special form of solution that helps us obtain the optimal approximate
	measure via simple iterative updates.
	
	\begin{theorem}\label{theoremGeneral}
		Assume that the approximate probability measure in problem (\ref{mipro1}) satisfies Assumptions \ref{assumption1}
		and the assumptions presented in Theorem \ref{generalAiClosed}.
		Using the same notations as in Theorem \ref{generalAiClosed}, in addition, we assume
		\begin{align}\label{assMain}
		\sup_{x_i\in T_N^i}\sup_{\substack{\nu^j\in\mathcal{A}_j \\ j\neq i}}
		\int_{\prod_{j\neq i} \mathcal{H}_j} \left( \Phi^0(x) + \Phi(x) \right)  1_A(x)
		\prod_{j\neq i}\nu^{j}(dx_j) < \infty,
		\end{align}
		and
		\begin{align}\label{assMain1}
		\sup_{\substack{\nu^j\in\mathcal{A}_j \\ j\neq i}} \!\int_{\mathcal{H}_i}\!\exp\bigg(\!\! -\int_{\prod_{j\neq i}\mathcal{H}_j}\!\!\!(\Phi^0(x)+\Phi(x))1_{A^c}(x)\prod_{j\neq i}\nu^{j}(dx_j)\!\bigg)M_i(x)\mu_{r}^{i}(dx_i) < \infty,
		\end{align}
		where $A:=\{ x \,|\, \Phi^0(x)+\Phi(x) \geq 0 \}$,
		and $M_i(x):=\max\left(1, a_i(\epsilon, x_i)\right)$ with $i,j=1,2,\cdots,M$.
		Then, problem (\ref{mipro1}) possesses a solution $\nu = \prod_{j=1}^{M}\nu^{j}\in\mathcal{M}(\mathcal{H})$ with the following form
		\mye{
			\frac{d\nu}{d\mu_{r}} \propto \exp\bigg( -\sum_{i=1}^{M}\Phi_{i}^{r}(x_{i}) \bigg),
		}
		where
		\mye{\label{mainThSpec1}
			\Phi_{i}^{r}(x_{i}) & = \int_{\prod_{j\neq i}\mathcal{H}_{j}}\!\! \bigg(\Phi^{0}(x) +
			\Phi(x) \bigg) \prod_{j\neq i}\nu^{j}(dx_{j}) + \text{Const}
		}
		and
		\mye{
			\nu^{i}(dx_{i}) & \propto \exp\big( -\Phi_{i}^{r}(x_{i}) \big)\mu_{r}^{i}(dx_{i}).
		}
	\end{theorem}
	
	\begin{remark}\label{conditionExplain1}
		For $i=1,2,\cdots, M$, conditions (\ref{assMain}) and (\ref{assMain1}) ensure that each components of
		the approximate measure $\nu$ and the reference probability measure $\mu_r$ are equivalent.
		These two conditions can be verified straightforwardly for specific examples relying on the
		integrability and boundedness conditions of $\Phi_i^r$ contained in the definition of $\mathcal{A}_i$ in (\ref{defAj})
		for $i=1,2,\cdots,M$.
	\end{remark}
	
	\begin{remark}\label{mainIteTheRemark}
		Formula (\ref{mainThSpec1}) means that the logarithm of the optimal solution for factor $\nu^{j}$
		can be obtained simply by considering the logarithm of the joint distribution over all of the other variables and then
		taking the expectation with respect to all of the other factors $\{ \nu^{i} \}$ fixed for $i\neq j$.
		This result is in accordance with the finite-dimensional case illustrated in Subsection 2.3 of \cite{Zhang2018IEEE}.
	\end{remark}
	
	\begin{remark}\label{algCon1Remark}
		Based on Theorem \ref{theoremGeneral}, we can therefore seek a solution by first initializing all
		of the potentials $\Phi_{j}^{r}$ appropriately and then cycling through the potentials and replacing each in turn
		with a revised estimate given by the right-hand side of (\ref{mainThSpec1}) evaluated by using the current
		estimates for all of the other potentials.
	\end{remark}
	
	\section{Applications to some general inverse problems}\label{SectionTypical}
	
	In Subsection \ref{SectionLinearIPG}, we apply our general theory to an abstract linear inverse problem (ALIP).
	We assume that the prior and noise probability measures are all Gaussian with some hyper-parameters,
	and then we formulate hierarchical models that can be efficiently solved by using the variational Bayes' approach.
	In Subsection \ref{SectionLinearIPL}, we assume that the noise is distributed according to the Laplace distribution.
	Through this assumption, we can formulate algorithms that solve ALIP and are robust to outliers.
	
	\subsection{Linear inverse problems with Gaussian noise}\label{SectionLinearIPG}
	In this subsection, we apply our general theory to an abstract linear inverse problem.
	A detailed investigation of the corresponding finite-dimensional case can be found in \cite{Jin2010JCP}.
	
	Let $\mathcal{H}_{u}$ be some separable Hilbert space and $N_{d}$ be a positive integer.
	We describe the linear inverse problem as
	\mye{\label{linearAIP}
		d = Hu + \epsilon,
	}
	where $d\in\mathbb{R}^{N_{d}}$ is the measurement data, $u\in \mathcal{H}_{u}$ is the sought-for solution,
	$H$ is a bounded linear operator from $\mathcal{H}_{u}$ to $\mathbb{R}^{N_{d}}$,
	and $\epsilon$ is a Gaussian random vector with zero mean and $\tau^{-1}I$ variance.
	We will focus on the hyper-parameter treatment within hierarchical models and the challenges in efficiently exploring
	the posterior probability.
	
	To formulate this problem under the Bayesian inverse framework, we introduce a prior probability
	measure for the unknown function $u$. Let $\mathcal{C}_{0}$ be a symmetric, positive definite and trace class
	operator defined on $\mathcal{H}_{u}$, and let $(e_{k}, \alpha_{k})$ be an eigen-system of the operator $\mathcal{C}_{0}$ such that
	$\mathcal{C}_{0}e_{k} = \alpha_{k}e_{k}$.
	Without loss of generality, we assume that the eigenvectors $\{e_k\}_{k=1}^{\infty}$ are orthonormal and the
	eigenvalues $\{\alpha_k\}_{k=1}^{\infty}$ are in a descending order.
	In the following, for a function $u\in\mathcal{H}_u$, we denote $u_j := \langle u,e_j \rangle$ for $j=1,2,\cdots$. 
	According to Subsection 2.4 in \cite{Dashti2014}, we have
	\mye{
		\mathcal{C}_{0} = \sum_{j = 1}^{\infty}\alpha_{j}e_{j}\otimes e_{j},
	}
	where $\otimes$ denotes the tensor product on Hilbert space \cite{Lax2002FA,Reed1980FunctionalAnalysis}.
	As indicated in \cite{Cotter2013SS,Cui2016JCP,Feng2018SIAM}, we assume that the data are only
	informative on a finite number of directions in $\mathcal{H}_{u}$.
	Under this assumption, we introduce a positive integer $K$, which represents the number of dimensions
	that is informed by the data (i.e., the so-called intrinsic dimensionality), which is different
	from the discretization dimensionality, i.e., the number of mesh points used to represent the unknown variables.
	The value of $K$ can be specified with a heuristic approach \cite{Feng2018SIAM}:
	\begin{align}
	K = \min\left\{ k\in\mathbb{N} \,\bigg|\, \frac{\alpha_k}{\alpha_1}<\epsilon \right\},
	\end{align}
	where $\epsilon$ is a prescribed threshold.
	Let $\lambda$ be a positive real number, then, we define
	\mye{
		\mathcal{C}_{0}^{K}(\lambda) := \sum_{j=1}^{K}\lambda^{-1}\alpha_{j}e_{j}\otimes e_{j}
		+ \sum_{j=K+1}^{\infty}\alpha_{j}e_{j}\otimes e_{j},
	}
	which is obviously a symmetric, positive definite and trace-class operator.
	Numerical results shown in \cite{Feng2018SIAM} indicate that the above heuristic approach could provide
	acceptable results when $\epsilon$ is small enough for a lot of practical inverse problems.
	However, if the data is particularly informative and far from the prior,
	this heuristic approach may lead to a Bayesian inference model that has not
	enough ability to incorporate information encoded in data.
	Concerned with this problem, we intend to give more detailed discussions in our future work.
	We refer to some recent studies \cite{Agapiou2014SIAM,NIPS2019_9649} that may provide some useful ideas.
	Then, we assume
	\mye{
		u \sim \mu_{0}^{u,\lambda} = \mathcal{N}(u_{0}, \mathcal{C}_{0}^{K}(\lambda)).
	}
	Let $Gamma(\alpha, \beta)$ be the Gamma probability measure defined on $\mathbb{R}^{+}$ with
	the probability density function $p_{G}$ expressed as
	\mye{
		p_{G}(x;\alpha,\beta) = \frac{\beta^{\alpha}}{\Gamma(\alpha)}x^{\alpha-1}e^{-\beta x},
	}
	where $\Gamma(\cdot)$ is the usual Gamma function.
	Then, except for the function $u$, we assume that the parameters $\lambda$ and $\tau$ involved in the prior and
	noise probability measures are all random variables satisfying
	$\lambda \sim \mu_{0}^{\lambda} := Gamma(\alpha_{0}, \beta_{0})$ and
	$\tau \sim \mu_{0}^{\tau} := Gamma(\alpha_{1}, \beta_{1}).$
	With these preparations, we define the prior probability measure employed for this problem as follows:
	\mye{\label{Ex1PriorDef}
		\mu_{0}(du, d\lambda, d\tau) = \mu_{0}^{u,\lambda}(du)\mu_{0}^{\lambda}(d\lambda)\mu_{0}^{\tau}(d\tau).
	}
	Let $\mu$ be the posterior probability measure for random variables $u$, $\lambda$, and $\tau$.
	According to Theorems 15 and 16 proved in \cite{Dashti2014}, this probability measure can be defined as
	\mye{\label{beijia1}
		\frac{d\mu}{d\mu_{0}}(u,\lambda,\tau) = \frac{1}{Z_{\mu}} \tau^{N_{d}/2} \exp\bigg( -\frac{\tau}{2}\|Hu-d\|^{2} \bigg),
	}
	where
	\mye{\label{formZmu1}
		Z_{\mu} = \int_{\mathcal{H}_{u}\times \mathbb{R}^{+}\times \mathbb{R}^{+}}
		\tau^{N_{d}/2}\exp\bigg( -\frac{\tau}{2}\|Hu-d\|^{2} \bigg)\mu_{0}(du, d\lambda, d\tau).
	}
	
	To apply the general theory developed in Section \ref{SectionGeneralTheory}, we specify the following
	reference probability measure
	\footnote{In practical machine learning applications, especially for large-scale scenarios, researchers often assume the approximating measures are independent in each component (a fully diagonal approximation to the posterior) that further reduce of computational burden. This, however, also tends to decrease the computation accuracy due to the neglecting of existing correlations between different components of $u$. We thus consider to preserve such correlation in our method to alleviate the possible negative influence of ignoring such beneficial knowledge.}
	\mye{
		\mu_{r}(du, d\lambda, d\tau) = \mu_{r}^{u}(du) \mu_{r}^{\lambda}(d\lambda) \mu_{r}^{\tau}(d\tau),
	}
	where $\mu_{r}^{u}=\mathcal{N}(u_{0}, \mathcal{C}_{0})$ is a Gaussian probability measure,
	and $\mu_{r}^{\lambda}$ and $\mu_{r}^{\tau}$ are chosen to be $\mu_{0}^{\lambda}$ and $\mu_{0}^{\tau}$, respectively.
	
	In Assumption \ref{assumption1}, we assume that the approximate probability measure is separable
	with respect to the random variables $u$, $\lambda$, and $\tau$ with the form
	\mye{\label{Ex1NU1}
		\nu(du,d\lambda,d\tau) = \nu^{u}(du)\nu^{\lambda}(d\lambda)\nu^{\tau}(d\tau).
	}
	In addition, we assume that its Radon-Nikodym derivative with respect to $\mu_{r}$ can be written as
	\mye{
		\frac{d\nu}{d\mu_{r}}(u, \lambda, \tau) = \frac{1}{Z_{r}}
		\exp\bigg( -\Phi_{u}^{r}(u)-\Phi_{\lambda}^{r}(\lambda)-\Phi_{\tau}^{r}(\tau) \bigg).
	}
	For the Radon-Nikodym derivative of $\mu_{0}$ with respect to $\mu_{r}$, we have
	\mye{\label{linerEx1}
		\frac{d\mu_{0}}{d\mu_{r}}(u,\lambda,\tau) & = \frac{d\mu_{0}^{u,\lambda}}{d\mu_{r}^{u}}(u)
		\frac{d\mu_{0}^{\lambda}}{d\mu_{r}^{\lambda}}(\lambda)\frac{d\mu_{0}^{\tau}}{d\mu_{r}^{\tau}}(\tau) \\
		& = \lambda^{K/2}
		\exp\bigg( -\frac{1}{2}\|(\mathcal{C}_{0}^{K}(\lambda))^{-1/2}(u-u_0)\|^{2}
		+ \frac{1}{2}\|\mathcal{C}_{0}^{-1/2}(u-u_0)\|^{2} \bigg) \\
		& = \lambda^{K/2}\exp\bigg( -\frac{1}{2}\sum_{j=1}^{K}(u_{j}-u_{0j})^{2}(\lambda-1)\alpha_{j}^{-1} \bigg),
	}
	which implies that $\Phi^{0}$ introduced in Assumption \ref{assumption1} takes the following form:
	\mye{\label{jia1}
		\Phi^{0}(u,\lambda,\tau) =\frac{1}{2}\sum_{j=1}^{K}(u_{j}-u_{0j})^{2}(\lambda-1)\alpha_{j}^{-1} - \frac{K}{2}\log\lambda.
	}
	
	\begin{remark}\label{remarkEx1}
		It should be noted that $\mathbb{R}^{+}$ is not a Hilbert space. However, the general theory is constructed on some separable Hilbert spaces.
		This issue can be resolved by considering $\lambda' := \log\lambda$ and $\tau':= \log\tau$ instead of $\lambda$ and $\tau$.
		Through this simple transformation, the space of hyper-parameters becomes $\mathbb{R}$ which
		is a Hilbert space. The calculations presented here also hold true when considering $\lambda'$ and $\tau'$ as
		hyper-parameters. Actually, we can derive that $e^{\lambda'}$ and $e^{\tau'}$ are distributed according to the same Gamma distributions
		as $\lambda$ and $\tau$.
		Choosing $a_u(\epsilon, u), a_{\lambda'}(\epsilon,\lambda')$, and $a_{\tau'}(\epsilon,\tau')$ appropriately,
		we can verify the conditions proposed in Theorem \ref{theoremGeneral} 
		(critical steps are provided in the supplementary materials).
		In the following, we still use $\lambda$ and $\tau$ as hyper-parameters.
		With this a little abusive use of the general theory (can be rigorously verified through the above simple transformation),
		the reader may see more clearly the connections between the finite- and infinite-dimensional theory.
	\end{remark}
	
	We now calculate $\Phi_{u}^{r}(u)$, $\Phi_{\lambda}^{r}(\lambda)$, and $\Phi_{\tau}^{r}(\tau)$
	according to the general results as shown in Theorem \ref{theoremGeneral}.
	
	\textbf{Calculate $\Phi_{u}^{r}(u)$}:
	A direct application of formula (\ref{mainThSpec1}) yields
	\mye{\label{formulaPhiur}
		\Phi_{u}^{r}(u) = & \int_{0}^{\infty}\!\!\!\int_{0}^{\infty}\bigg(\frac{1}{2}
		\sum_{j=1}^{K}(u_{j}-u_{0j})^{2}(\lambda-1)\alpha_{j}^{-1} + \frac{\tau}{2}\|Hu-d\|^{2} \\
		& \qquad\qquad\qquad
		- \frac{K}{2}\log\lambda - \frac{N_d}{2}\log\tau \bigg)\nu^{\tau}(d\tau)\nu^{\lambda}(d\lambda) + \text{Const}	\\
		= & \frac{1}{2}\tau^{*}\|Hu-d\|^{2} + \frac{1}{2}(\lambda^{*}-1)\sum_{j=1}^{K}\alpha_{j}^{-1}(u_{j}-u_{0j})^{2} + \text{Const},
	}
	where
	\mye{\label{taulambdastar1}
		\tau^{*} = \mathbb{E}^{\nu^{\tau}}[\tau] = \int_{0}^{\infty}\tau\nu^{\tau}(d\tau) \quad \text{and} \quad
		\lambda^{*} = \mathbb{E}^{\nu^{\lambda}}[\lambda] = \int_{0}^{\infty}\lambda\nu^{\lambda}(d\lambda).
	}
	On the basis of equality (\ref{formulaPhiur}), we derive
	\mye{
		\frac{d\nu^{u}}{d\mu_{r}^{u}}(u) \propto \exp\bigg(
		-\frac{\tau^{*}}{2}\|Hu-d\|^{2} - \frac{\lambda^{*}-1}{2}\sum_{j=1}^{K}\alpha_{j}^{-1}(u_{j}-u_{0j})^{2} \bigg).
	}
	We define
	\mye{
		\mathcal{C}_{0}(\lambda^{*}) = \sum_{j=1}^{K}(\lambda^{*})^{-1}\alpha_{j}e_{j}\otimes e_{j}
		+ \sum_{j=K+1}^{\infty}\alpha_{j}e_{j}\otimes e_{j}.
	}
	Then, according to Example 6.23 in \cite{Stuart2010AN}, we know that the probability measure $\nu^{u}$
	is a Gaussian measure $\mathcal{N}(u^{*}, \mathcal{C})$ with
	\mye{\label{Ex1Formula1}
		\mathcal{C}^{-1}=\tau^{*}H^{*}H + \mathcal{C}_{0}(\lambda^{*})^{-1}
		\quad \text{and} \quad u^{*} = \mathcal{C}\big(\tau^{*}H^{*}d + \mathcal{C}_{0}(\lambda^{*})^{-1}u_{0}\big).
	}
	
	\textbf{Calculate $\Phi_{\lambda}^{r}(\lambda)$ and $\Phi_{\tau}^{r}(\tau)$}:
	According to formula (\ref{mainThSpec1}), we have
	\mye{
		\Phi_{\lambda}^{r}(\lambda) = & \int_{0}^{\infty}\!\!\!\int_{\mathcal{H}_{u}} \bigg(
		\frac{1}{2}\sum_{j=1}^{K}(u_{j}-u_{0j})^{2}\alpha_{j}^{-1}\lambda - \frac{K}{2}\log\lambda
		\bigg)\nu^{u}(du)\nu^{\tau}(d\tau) + \text{Const} \\
		= & \frac{1}{2}\mathbb{E}^{\nu^{u}}\bigg(\sum_{j=1}^{K}(u_{j}-u_{0j})^{2}\alpha_{j}^{-1}\bigg)\lambda
		- \frac{K}{2}\log\lambda + \text{Const},
	}
	which implies that
	\mye{
		\frac{d\nu^{\lambda}}{d\mu_{r}^{\lambda}}(\lambda) \propto
		\lambda^{K/2} \exp\bigg( -\frac{1}{2}\mathbb{E}^{\nu^{u}}
		\bigg( \sum_{j=1}^{K}(u_{j}-u_{0j})^{2}\alpha_{j}^{-1} \bigg)\lambda \bigg).
	}
	Given that $\lambda$ is a scalar random variable, we can write the density function as bellow:
	\mye{\label{Ex1Formula2}
		\rho_{G}(\lambda; \tilde{\alpha}_{0}, \tilde{\beta}_{0}) =
		\frac{\tilde{\beta}_{0}^{\tilde{\alpha}_{0}}}{\Gamma(\tilde{\alpha}_{0})}\lambda^{\tilde{\alpha}_{0}-1}
		\exp(-\tilde{\beta}_{0}\lambda),
	}
	where
	\mye{\label{Ex1Formula3}
		\tilde{\alpha}_{0} = \alpha_{0}+\frac{K}{2} \quad \text{and} \quad
		\tilde{\beta}_{0} = \beta_{0} + \frac{1}{2}\mathbb{E}^{\nu^{u}}
		\bigg( \sum_{j=1}^{K}(u_{j}-u_{0j})^{2}\alpha_{j}^{-1} \bigg).
	}
	Similar to the above calculations of $\Phi_{\lambda}^{r}(\lambda)$, we derive
	\mye{
		\Phi_{\tau}^{r}(\tau) = & \int_{0}^{\infty}\!\!\!\int_{\mathcal{H}_{u}}\bigg(
		\frac{\tau}{2}\|Hu-d\|^{2} - \frac{N_{d}}{2}\log\tau \bigg)\nu^{u}(du)\nu^{\lambda}(d\lambda) + \text{Const} \\
		= & \frac{1}{2}\mathbb{E}^{\nu^{u}}(\|Hu-d\|^{2})\tau - \frac{N_{d}}{2}\log\tau + \text{Const},
	}
	which implies
	\mye{
		\frac{d\nu^{\tau}}{d\mu_{r}^{\tau}}(\tau) \propto \tau^{\frac{N_{d}}{2}}
		\exp\bigg( -\frac{1}{2}\mathbb{E}^{\nu^{u}}(\|Hu-d\|^{2})\tau \bigg).
	}
	Therefore, $\nu^{\tau}$ is a Gamma distribution $Gamma(\tilde{\alpha}_{1}, \tilde{\beta}_{1})$ with
	\mye{\label{Ex1Formula4}
		\tilde{\alpha}_{1} = \alpha_{1} + \frac{N_{d}}{2} \quad \text{and} \quad
		\tilde{\beta}_{1} = \beta_{1} + \frac{1}{2}\mathbb{E}^{\nu^{u}}(\|Hu-d\|^{2}).
	}
	
	According to the statements in Remark \ref{algCon1Remark},
	we provide an iterative algorithm based on formulas (\ref{Ex1Formula1}), (\ref{Ex1Formula2}),
	(\ref{Ex1Formula3}), and (\ref{Ex1Formula4}) in Algorithm \ref{algComplexVBG}.
	\begin{algorithm}
		\setstretch{1}
		\caption{Variational approximation for the case of Gaussian noise}
		\label{algComplexVBG}
		\begin{algorithmic}
			\STATE {1: Give an initial guess $\mu_{0}^{u,\lambda}$ ($u_0$ and $\lambda$),
				$\mu_{0}^{\lambda}$ ($\alpha_0$ and $\beta_0$) and $\mu_{0}^{\tau}$ ($\alpha_1$ and $\beta_1$).\\
				\quad\,\! Specify the tolerance \emph{tol} and set $k=1$.
			}
			\STATE {2: \textbf{repeat}
			}
			\STATE {3: \qquad Set $k = k + 1$
			}
			\STATE {3: \qquad Calculate $\lambda_{k} = \mathbb{E}^{\nu^{\lambda}_{k-1}}[\lambda]$,
				$\tau_{k} = \mathbb{E}^{\nu^{\tau}_{k-1}}[\tau]$
			}
			\STATE {4: \qquad Calculate $\nu_{k}^{u}$ by
				\vskip -0.5 cm
				\begin{align*}
				\mathcal{C}_{k}^{-1} = \tau_{k} H^{*}H + \mathcal{C}_{0}(\lambda_{k})^{-1}, \quad
				u_{k} = \mathcal{C}_{k}\big( \tau_{k}H^{*}d + \mathcal{C}_{0}(\lambda_{k})^{-1}u_{0} \big).
				\end{align*}
			}\vskip -1 cm
			\STATE{5: \qquad Calculate $\nu_{k}^{\lambda}$ and $\nu_{k}^{\tau}$ by
				\vskip -0.5 cm
				\myen{
					\nu_{k}^{\lambda} = Gamma(\tilde{\alpha}_{0}, \tilde{\beta}_{0}^{k}), \quad
					\nu_{k}^{\tau} = Gamma(\tilde{\alpha}_{1}, \tilde{\beta}_{1}^{k}),
				}
				\vskip -0.2 cm
				\qquad\quad\,\! where
				\vskip -0.8 cm
				\myen{
					&\tilde{\alpha}_{0} = \alpha_0 + \frac{K}{2}, \quad
					\tilde{\beta}_{0}^{k} = \beta_{0} + \frac{1}{2}\mathbb{E}^{\nu_{k}^{u}}
					\bigg( \sum_{j=1}^{K}(u_{j}-u_{0j})^{2}\alpha_{j}^{-1} \bigg), \\
					& \tilde{\alpha}_{1} = \alpha_{1} + \frac{N_{d}}{2}, \quad
					\tilde{\beta}_{1}^{k} = \beta_{1} + \frac{1}{2}\mathbb{E}^{\nu_{k}^{u}}(\|Hu-d\|^{2}).
				}\vskip -0.5 cm }
			\STATE{6: \textbf{until} $\max\big(\|u_k - u_{k-1}\|/\|u_k\|, \|\lambda_k - \lambda_{k-1}\|/\|\lambda_k\|,
				\|\tau_k - \tau_{k-1}\|/\|\tau_k\| \big) \leq$ \emph{tol}
			}
			\STATE{7: Return $\nu_{k}^{u}(du)\nu_{k}^{\lambda}(d\lambda)\nu_{k}^{\tau}(d\tau)$
				as the solution.
			}
		\end{algorithmic}
	\end{algorithm}
	Before giving the next example, we provide a brief discussion of the computational details and the cost of Algorithm 1.
	For small- or medium-scale problems, we may construct the finite-dimensional approximate operators $H$ and $H^*$ explicitly \cite{Jin2010JCP}.
	However, for large-scale problems, it is impossible to build finite-dimensional approximations explicitly.
	Actually, for running the iterations, we only need to compute the mean estimates $u_k$ and some quantities related to $\nu_k^u$ such as
	$\mathbb{E}^{\nu_k^u}(\|Hu-d\|^2)$. For obtaining mean estimates, we can use
	a matrix-free conjugate gradient (CG) method \cite{Tan2013IP,Golub1996book,Petra2014SIAM} to solve the following problem
	\begin{align}\label{conj1}
	(\tau_k H^* H + \mathcal{C}_0(\lambda_k)^{-1})u_k = \tau_k H^* d + \mathcal{C}_0(\lambda_k)^{-1}u_0,
	\end{align}
	where no explicit forms of $H^*H$ and $H^*$ need to be constructed.
	As demonstrated in \cite{Tan2013IP,Petra2014SIAM}, the CG iterations may
	be terminated when sufficient reduction is made in the norm of the gradient and the prior operator may also be used to precondition 
	the CG iterations.
	For the term $\mathbb{E}^{\nu_k^u}(\|Hu-d\|^2)$, by a straightforward generalization of the finite dimensional case \cite{Jin2010JCP}
	(Proposition 1.18 in \cite{Prato2006IDAnalysis} and (c) of Theorem VI.25 in \cite{Reed1980FunctionalAnalysis} are used),
	we know that the core difficulty is to compute the following quantity
	\begin{align}\label{trance1}
	\text{Tr}\big((\tau_k \mathcal{C}_0(\lambda_k)^{1/2}H^* H\mathcal{C}_0(\lambda_k)^{1/2} + Id)^{-1}
	\mathcal{C}_0(\lambda_k)^{1/2}H^* H\mathcal{C}_0(\lambda_k)^{1/2}\big),
	\end{align}
	where $\text{Tr}(\cdot)$ represents taking trace of the operator.
	For a lot of practical applications, the operator $H^*H$ is a compact operator.
	Then the analysis provided in Subsections 5.2 and 5.4 in \cite{Tan2013IP} may applicable in the current setting,
	which implies that only a small number of eignvalues (independent of the dimension of the discretized parameter field)
	is required to be evaluated. We intend to discuss in detail on efficient implementations for large-scale problems in our future work.
	
	\subsection{Linear inverse problems with Laplace noise}\label{SectionLinearIPL}
	As revealed by previous studies on low-rank matrix factorization \cite{Zhao2015IEEE}, the
	Gaussian noise tends to be sensitive to outliers. Compared with the Gaussian distribution, the Laplace
	distribution is a \emph{heavy-tailed} distribution that can better fit heavy noises and outliers.
	In this subsection, we develop VBM for the linear
	inverse problem (\ref{linearAIP}) with the Laplace noise assumption.
	
	For the noise vector
	$\epsilon = (\epsilon_1, \epsilon_2, \cdots, \epsilon_{N_d})^{T} \in \mathbb{R}^{N_{d}},$
	we assume that each component $\epsilon_{i}$ follows the Laplace distribution with zero mean
	\mye{\label{noisedis1}
		\epsilon_{i} \sim \text{Laplace}\bigg(0, \sqrt{\frac{\tau}{2}}\bigg)
	}
	with $\tau\in\mathbb{R}^{+}$. The probability density function of the above Laplace distribution
	is denoted by $p_{L}(\epsilon_{i}|0, \sqrt{\tau/2})$ that takes the following form:
	\mye{
		p_{L}(\epsilon_{i} | 0,\sqrt{\tau/2})  = \sqrt{\frac{2}{\tau}}
		\exp\bigg( -\frac{|\epsilon_{i}|}{\sqrt{\tau /2}} \bigg) .
	}
	
	However, the Laplace distribution cannot be easily employed for posterior inference
	within the variational Bayes' inference framework \cite{Zhao2015IEEE}.
	A commonly utilized strategy will be employed to reformulate the Laplace distribution as a Gaussian scale mixture with
	exponential distributed prior to the variance, as discussed in \cite{Andrews1974JRSSSB,Zhao2015IEEE}.
	Let $p_{E}(z|\tau)$ be the density function of an exponential distribution, that is,
	\mye{
		p_{E}(z|\tau) = \frac{1}{\tau}\exp\bigg( -\frac{z}{\tau} \bigg).
	}
	Then, we have
	\mye{
		p_{L}\bigg(x\big|0, \sqrt{\frac{\tau}{2}}\bigg) = & \frac{1}{2}\sqrt{\frac{2}{\tau}}
		\exp\bigg( -\sqrt{\frac{2}{\tau}}|x| \bigg)	\\
		= & \int_{0}^{\infty}\frac{1}{\sqrt{2\pi z}}\exp\bigg( -\frac{x^2}{2z} \bigg)\frac{1}{\tau}
		\exp\bigg( -\frac{z}{\tau} \bigg)dz \\
		= & \int_{0}^{\infty}p_{N}(x|0,z)p_{E}(z|\tau)dz.
	}
	By substituting (\ref{noisedis1}) into the above equation, we obtain
	\mye{
		p_{L}(\epsilon_{i} | 0, \sqrt{\tau/2}) = \int_{0}^{\infty}p_{N}(\epsilon_{i}|0,z_{i})p_{E}(z_{i}|\tau)dz_{i},
	}
	where $p_{N}(\epsilon_{i}|0,z_{i})$ is the density function of a Gaussian measure on $\mathbb{R}$ with
	a zero mean and $z_{i}$ variance.
	Thus, we can impose a two-level hierarchical prior instead of a single-level Laplace prior on each $\epsilon_{i}$ as
	\mye{\label{NoiseModel1}
		\epsilon_{i}\sim\mathcal{N}(0,z_{i}), \qquad z_{i}\sim\text{Exponential}(\tau).
	}
	Let $w_{i} = z_{i}^{-1}$. Given that $z_{i} \sim \text{Exponential}(\tau)$, we know that $w_{i}\sim\mu_{0}^{w_{i}}$ with
	$\mu_{0}^{w_{i}}$ being a probability distribution with the following probability density function:
	\mye{
		\frac{1}{\tau}\exp\bigg(-\frac{1}{\tau w_{i}}\bigg)\frac{1}{w_{i}^{2}}.
	}
	Let $W$ be a diagonal matrix with diagonal $w = \{ w_{1}, w_{2}, \cdots, w_{N_{d}} \}$, and let
	\mye{\label{priorLap2}
		\mu_{0}^{w} = \prod_{i=1}^{N_{d}}\mu_{0}^{w_{i}}.
	}
	For the prior probability measure of $u$, similar to Subsection \ref{SectionLinearIPG},
	we set this measure as for the Gaussian noise case, that is,
	\mye{\label{priorLap1}
		u \sim \mu_{0}^{u,\lambda}=\mathcal{N}(u_0, \mathcal{C}_{0}^{K}(\lambda)), \quad
		\lambda \sim \mu_{0}^{\lambda} = \text{Gamma}(\alpha_{0}, \beta_{0}).
	}
	By combining (\ref{priorLap2}) and (\ref{priorLap1}), we obtain the full prior probability measure as
	\mye{
		\mu_{0}(du, d\lambda, dw) = \mu_{0}^{u,\lambda}(du)\mu_{0}^{\lambda}(d\lambda)\mu_{0}^{w}(dw).
	}
	For the reference probability measure, we set
	$\mu_{r}(du, d\lambda, dw) = \mu_{r}^{u}(du)\mu_{r}^{\lambda}(d\lambda)\mu_{r}^{w}(dw)$,
	where $\mu_{r}^{u}=\mathcal{N}(u_0,\mathcal{C}_{0})$, $\mu_{r}^{\lambda} = \mu_{0}^{\lambda}$, and
	$\mu_{r}^{w}=\mu_{0}^{w}$. By similar calculations as shown in (\ref{linerEx1}), we obtain
	\mye{
		\Phi^{0}(u,\lambda,\tau) =\frac{1}{2}\sum_{j=1}^{K}(u_{j}-u_{0j})^{2}(\lambda-1)\alpha_{j}^{-1} - \frac{K}{2}\log\lambda.
	}
	For the posterior probability measure, by assumptions on the noises (\ref{NoiseModel1})-(\ref{priorLap2}), we have
	\mye{
		\frac{d\mu}{d\mu_{0}}(u, \lambda, w) = \frac{1}{Z_{\mu}}|W|^{1/2}
		\exp\bigg( -\frac{1}{2}\|W^{1/2}(Hu-d)\|^{2} \bigg),
	}
	which implies
	$\Phi(u,\lambda,w) = \frac{1}{2}\|W^{1/2}(Hu-d)\|^2 - \frac{1}{2}\log|W|$.
	Similar to the Gaussian noise case, we specify the approximate probability measure as
	\mye{\label{Ex2AppMea}
		\frac{d\nu}{d\mu_{r}}(u, \lambda, w) = \frac{1}{Z_{r}}
		\exp\bigg( -\Phi_{u}^{r}(u)-\Phi_{\lambda}^{r}(\lambda)-\Phi_{w}^{r}(w) \bigg).
	}
	
	With these preparations, we are ready to calculate the three potentials in (\ref{Ex2AppMea}).
	As discussed in Remark \ref{remarkEx1}, we use $\lambda > 0$ as a hyper-parameter, which is not in accordance
	with our general theory. However, it can be made rigorous by considering $\lambda'=\ln\lambda$ as the hyper-parameter.
	
	\textbf{Calculate $\Phi_{u}^{r}$}:
	Following formula (\ref{mainThSpec1}), we can derive
	\mye{\label{LapCal1}
		\Phi_{u}^{r}(u) = & \int\!\!\!\!\int
		\frac{1}{2}\sum_{j=1}^{K}(u_{j}-u_{0j})^{2}(\lambda-1)\alpha_{j}^{-1}\! + \!\frac{1}{2}\|W^{1/2}(Hu-d)\|^{2}
		d\nu^{\lambda}d\nu^{w} \! + \! \text{Const}	\\
		= & \frac{\lambda^{*}-1}{2}\sum_{j=1}^{K}\alpha_{j}^{-1}(u_{j}-u_{0j})^{2}
		+ \frac{1}{2}\|W^{*}(Hu-d)\|^{2} + \text{Const},
	}
	where
	$\lambda^{*} = \mathbb{E}^{\nu^{\lambda}}[\lambda]$ and
	$W^{*} = \text{diag}(\mathbb{E}^{\nu^{w}}[w_1], \mathbb{E}^{\nu^{w}}[w_2], \cdots, \mathbb{E}^{\nu^{w}}[w_{N_{d}}])$.
	From the equality (\ref{LapCal1}), we easily conclude that
	\mye{
		\frac{d\nu^{u}}{d\mu_{r}^{u}}(u)\propto\exp\bigg(-\frac{1}{2}\|(W^{*})^{1/2}(Hu-d)\|^{2}
		- \frac{\lambda^{*}-1}{2} \sum_{j=1}^{K}\alpha_{j}^{-1}(u_{j}-u_{0j})^{2}  \bigg),
	}
	which implies that $u$ is distributed according to a Gaussian measure with a covariance operator and a mean value specified as
	$\mathcal{C}^{-1} = H^{*}W^{*}H + \mathcal{C}_{0}(\lambda^{*})^{-1}$ and
	$u^{*} = \mathcal{C}\big( H^{*}W^{*}d + \mathcal{C}_{0}(\lambda^{*})^{-1}u_{0} \big)$.
	
	\textbf{Calculate $\Phi_{\lambda}^{r}$}:
	Following formula (\ref{mainThSpec1}), we can derive
	\mye{
		\Phi_{\lambda}^{r}(\lambda) = & \int\!\!\!\int \frac{1}{2}\sum_{j=1}^{K}(u_{j}-u_{0j})^{2}\alpha_{j}^{-1}\lambda
		- \frac{K}{2}\log\lambda d\nu^{u}d\nu^{w} + \text{Const} \\
		= & \frac{1}{2}\mathbb{E}_{\nu^{u}}\bigg( \sum_{j=1}^{K}(u_{j}-u_{0j})^{2}\alpha_{j}^{-1} \bigg) \lambda
		- \frac{K}{2}\log\lambda + \text{Const}.
	}
	Therefore, we have
	\mye{
		\frac{d\nu^{\lambda}}{d\mu_{r}^{\lambda}}(\lambda) \propto \lambda^{K/2}\exp\bigg(
		-\frac{1}{2}\mathbb{E}^{\nu^{u}}\bigg( \sum_{j=1}^{K}(u_{j}-u_{0j})^{2}\alpha_{j}^{-1} \bigg) \lambda \bigg),
	}
	which implies that $\nu^{\lambda}$ is a Gamma distribution denoted by
	$\text{Gamma}(\tilde{\alpha}_{0}, \tilde{\beta}_{0})$ with
	\mye{
		\tilde{\alpha}_{0} = \alpha_{0}+K/2, \quad
		\tilde{\beta}_{0} = \beta_{0} + \frac{1}{2}\mathbb{E}^{\nu^{u}}\bigg( \sum_{j=1}^{K}
		(u_{j}-u_{0j})^{2}\alpha_{j}^{-1} \bigg).
	}
	
	\textbf{Calculate $\Phi_{w}^{r}$}:
	Following formula (\ref{mainThSpec1}), we derive
	\mye{
		\Phi_{w}^{r}(w)= & \int\!\!\!\int \frac{1}{2}\|W^{1/2}(Hu-d)\|^{2} - \frac{1}{2}\log|W| d\nu^{u}d\nu^{\lambda} + \text{Const} \\
		= & \frac{1}{2}\sum_{j=1}^{N_d}\mathbb{E}^{\nu^{u}}\big[ (Hu-d)_{i}^{2} \big] w_{i}
		- \frac{1}{2}\sum_{j=1}^{N_d}\log w_{i} + \text{Const},
	}
	which implies
	\mye{\label{LapCal2}
		\frac{d\nu^{w}}{d\mu_{r}^{w}}(w) \propto \prod_{j=1}^{N_d}w_{j}^{1/2}
		\exp\bigg( -\frac{1}{2}\mathbb{E}^{\nu^{u}}\big[ (Hu-d)_{j}^{2} \big] w_{j} \bigg).
	}
	Because $w$ is a finite dimensional random variable, we find
	\mye{
		d\nu^{w}\propto & \prod_{j=1}^{N_d}w_{j}^{1/2}
		\exp\bigg( -\frac{1}{2}\mathbb{E}^{\nu^{u}}\big[ (Hu-d)_{j}^{2} \big] w_{j} \bigg)
		\frac{1}{\tau}\exp\bigg( -\frac{1}{\tau w_{j}} \bigg) \frac{1}{w_{j}^{2}} dw \\
		\propto & \prod_{j=1}^{N_d} \frac{1}{\tau w_{j}^{3/2}}
		\exp\bigg( -\frac{1}{2}\mathbb{E}^{\nu^{u}}\big[ (Hu-d)_{j}^{2} \big] w_{j} -\frac{1}{\tau w_{j}} \bigg) dw.
	}
	In other words, $\nu^{w}$ is an inverse Gaussian distribution denoted by
	$\prod_{j=1}^{N_d}IG(m_{w_{j}}, \zeta)$ with
	\mye{
		m_{w_{j}} = \sqrt{\frac{2}{\tau \mathbb{E}^{\nu^{u}}\big[ (Hu-d)_{j}^{2} \big]}},
		\quad \zeta = \frac{2}{\tau}.
	}
	
	\textbf{Specify the parameter $\tau$}:
	From (\ref{NoiseModel1}), we know the parameter $\tau$ is directly related to noise variance parameter $z_{i} = w_{i}^{-1}$.
	Therefore, this parameter should be adjusted carefully to obtain reasonable results. Empirical Bayes \cite{Bishop2006PRML}
	provides an off-the-shelf tool to be adaptively tuned based on the noise information extracted from the data
	by updating it through
	$\tau = \frac{1}{N_d}\sum_{j=1}^{N_d}m_{w_j}^{-1} + \zeta^{-1}$.
	Using this elaborate tool, $\tau$ can be properly adapted to real data variance.
	
	Similar to the Gaussian noise case, an iterative algorithm, namely Algorithm 2, is constructed based on the above calculations.
	For large-scale problems, similar discussions of Algorithm 1 can be applied here.
	The only difference is on (\ref{trance1}) which is replaced by the following quantity
	\begin{align}\label{trance2}
	\text{Tr}\big((\tau_k \mathcal{C}_0(\lambda_k)^{1/2}H^* W_kH\mathcal{C}_0(\lambda_k)^{1/2} + Id)^{-1}
	\mathcal{C}_0(\lambda_k)^{1/2}H^*W_k H\mathcal{C}_0(\lambda_k)^{1/2}\big).
	\end{align}
	By our understanding, quantity (\ref{trance2}) could be calculated in a similar way as (\ref{trance1}). 
	\begin{algorithm}
		\caption{Variational approximation for the case of Laplace noise}
		\label{algComplexVBL}
		\begin{algorithmic}
			\STATE {1: Give an initial guess $\mu_{0}^{u,\lambda}$ ($u_0$ and $\lambda$),
				$\mu_{0}^{\lambda}$ ($\alpha_0$ and $\beta_0$), $\mu_{0}^{w}$ and $\tau$.\\
				\quad\,\! Specify the tolerance \emph{tol} and set $k=1$.
			}
			\STATE {2: \textbf{repeat}
			}
			\STATE {3: \qquad Set $k = k + 1$
			}
			\STATE {3: \qquad Calculate $\lambda_{k} = \mathbb{E}^{\nu^{\lambda}_{k-1}}[\lambda]$,
				$W_{k} = \text{diag}\big( \mathbb{E}^{\nu^{w}}[w_1], \mathbb{E}^{\nu^{w}}[w_2], \cdots,
				\mathbb{E}^{\nu^{w}}[w_{N_d}] \big)$ and \\
				\qquad\quad\, $\tau_{k} = \frac{1}{N_d}\sum_{j=1}^{N_d}(m_{w_j}^{k-1})^{-1} + (\zeta_{k-1})^{-1}$.
			}
			\STATE {4: \qquad Calculate $\nu_{k}^{u}$ by
				\myen{
					\mathcal{C}_{k}^{-1} = H^{*}W_{k}H + \mathcal{C}_{0}(\lambda_{k})^{-1}, \quad
					u_{k} = \mathcal{C}_{k}\big( H^{*}W_{k}d + \mathcal{C}_{0}(\lambda_{k})^{-1}u_{0} \big).
			}}
			\vspace{-0.5cm}
			\STATE{5: \qquad Calculate $\nu_{k}^{\lambda}$ and $\nu_{k}^{w}$ by
				\vspace{-0.5cm}
				\myen{
					& \nu_{k}^{\lambda} = Gamma(\tilde{\alpha}_{0}, \tilde{\beta}_{0}^{k}), \quad
					\nu_{k}^{w} = \prod_{j=1}^{N_{d}}IG(m_{w_{j}}^{k}, \zeta_{k}),  \\
					\tilde{\beta}_{0}^{k} = & \beta_{0} + \frac{1}{2}\mathbb{E}^{\nu_{k}^{u}}
					\bigg( \sum_{j=1}^{K}(u_{j}-u_{0j})^{2}\alpha_{j}^{-1} \bigg), \quad
					\tilde{\alpha}_{0} = \alpha_{0} + K/2, 	\\
					& \quad m_{w_{j}}^{k} = \sqrt{\frac{2}{\tau_{k} \mathbb{E}^{\nu_{k}^{u}}\big[ (Hu-d)_{j}^{2} \big]}},
					\quad \zeta_{k} = \frac{2}{\tau_{k}}.
				}
				\vskip -0.3 cm
			}
			\STATE{6: \textbf{until} $\max\big(\|u_k - u_{k-1}\|/\|u_k\|, \|\lambda_k - \lambda_{k-1}\|/\|\lambda_k\|,
				\|\tau_k - \tau_{k-1}\|/\|\tau_k\| \big) \leq$ \emph{tol}
			}
			\STATE{7: Return $\nu_{k}^{u}(du)\nu_{k}^{\lambda}(d\lambda)\nu_{k}^{w}(dw)$
				as the solution.
			}
		\end{algorithmic}
	\end{algorithm}

	%
	%
	
	\section{Concrete numerical examples}\label{SectionConcreteExamples}
	
	
	\subsection{Inverse source problem for Helmholtz equation}\label{SectionISPGeneral}
	
	The inverse source problem (ISP) studied in this section is borrowed from
	\cite{Bao2015IP,Bao2015SIAM,Cheng2016JDE,Isakov2018IPI},
	which determines the unknown current density function from measurements of the radiated fields at multiple wavenumbers.

	Consider the Helmholtz equation
	\mye{\label{HelmholtzEq}
		\Delta v + \kappa^{2}(1+q(x)) v = u_s\quad \text{in} \quad \mathbb{R}^{N_s},
	}
	where $N_s = 1,2$ is the space dimension, $\kappa$ is the wavenumber,
	$v$ is the radiated scalar field, and the source current density function $u_s(x)$
	is assumed to have a compact support.
	For the one-dimensional case, let the radiated field $v$ satisfy the
	absorbing boundary condition: $\partial_{r}v = i\kappa v$.
	For the two-dimensional case, let the radiated field $v$ satisfy the Sommerfeld radiation condition:
	$\partial_{r}v - i\kappa v = o(r^{-1/2}) \text{ as } r = |x|\rightarrow \infty$.
	In addition, we employ an uniaxial perfect match layer (PML) technique to truncate
	the whole plane into a bounded rectangular domain when $N_s = 2$.
	For details on the uniaxial PML technique, see \cite{Bao2010IP,Jia2019IP} and references therein.
	Let $D$ be the domain with absorbing layers, and $\Omega$ be the physical domain without absorbing layers.
	
	The ISP aims to determine the source function $u_s$ from the boundary measurements of the radiated field
	on the boundary $\partial\Omega$ for a series of wavenumbers. For clarity, we summarize the problem as follows:
	\begin{description}
		\item[Available data] For $0 < \kappa_{1} < \kappa_{2} < \cdots < \kappa_{N_f} < \infty$ ($N_f\in\mathbb{N}^{+}$), and
		measurement points $x^1,x^2,\cdots,x^{N_m} \in \partial\Omega$, we denote
		$$d^{\dag} := \big\{v(x^{i}, \kappa_{j}) \, | \, i = 1,2,\cdots,N_m, \text{ and } j=1,2,\cdots,N_f \big\}.$$
		The available data set is $d:=d^{\dag} + \epsilon$, where $\epsilon$ is the measurement error.
		\item[Unknown function] The source density function $u_s$ needs to be determined.
	\end{description}
	
	Generally, we let $\mathcal{F}_{\kappa}$ be the forward operator that maps $u_s$ to the solution $v$ when the wavenumber is $\kappa$ ,
	and let $\mathcal{M}$ be the measurement operator mapping $v$ to the available data.
	With these notations, the problem can be written abstractly as
	\mye{\label{HelmAbstract1}
		d_{\kappa} = H_{\kappa}(u_s) + \epsilon_{\kappa},
	}
	where $H_{\kappa}:=\mathcal{M}\circ\mathcal{F}_{\kappa}$ is the forward operator, and $\epsilon_{\kappa}$ is the random noise.
	
	
	To avoid inverse crime, we use a fine mesh to generate data and a rough mesh for the inversion.
	For the one-dimensional problem, meshes with mesh numbers of $1000$ and $600$ are
	used for the data generation and inversion, respectively.
	For the two-dimensional problem, we will provide details in the sequel.
	
	When the dimension of the parameters is relatively low, the proposed Algorithms 1 and 2
	are similar to the one build for the finite-dimensional case.
	Detailed comparisons with the MCMC algorithm have been given in \cite{Jin2012JCP,Jin2010JCP}, which reflect
	that high accurate inferences can be generated. Hence we will not present a comparison with the
	MCMC algorithm in the sequel for a relatively low dimensional case.
	For the infinite-dimensional Bayesian method with hyper-parameters,
	the noncentered algorithms are a more appropriate choice as illustrated in \cite{Agapiou2014SIAM}.
	Using the proposed general framework for the noncentered parameterize strategy and providing a comparison
	with the method proposed in \cite{Agapiou2014SIAM} could be an interesting future research problem. 
	
	It should be indicated that the finite element method is implemented by employing the open
	software FEniCS (Version 2018.1.0). For additional information on FEniCS, see \cite{Logg2012Book}.
	All programs were run on a personal computer with Intel(R) Core(TM) i7-7700 at 3.60 GHz (CPU), 32 GB (memory),
	and Ubuntu 18.04.2 LTS (OS).
	
	\subsection{One-dimensional ISP}\label{OneDISPSection}
	
	For clarity, we list the specific choices for some parameters introduced in Section \ref{SectionTypical} as follows:
	\begin{itemize}
		\item The operator $\mathcal{C}_{0}$ is chosen to be $(\text{Id}-\partial_{xx})^{-1}$ and taken $\epsilon=10^{-3}$.
		Here, the Laplace operator is defined on $\Omega$ with the zero Dirichlet boundary condition.
		\item The wavenumber series are specified as
		$\kappa_{j} = j \text{ with } j = \frac{1}{2}, 1, \frac{3}{2}, 2,\cdots,50$.
		\item Let domain $\Omega$ be an interval $[0,1]$, with $\partial\Omega = \{0, 1\}$. And the available data
		are assumed to be
		$\{ v(x^{i}, \kappa_{j}) \, | \, i = 1,2, \,\, x^1 = 0,\,\, x^2 = 1, \text{ and }j = 1,2,\cdots,100 \}$.
		\item The initial values required by Algorithm \ref{algComplexVBG} are chosen as
		$u_0 = 0, \alpha_{0} = \alpha_1 = 1, \beta_{0} = 10^{-1}, \beta_{1} = 10^{-5}$.
		The initial values required by Algorithm \ref{algComplexVBL} are chosen as
		$u_0 = 0, \alpha_{0} = 1, \beta_{0} = 10^{-1}, \tau = 10^{-7}$.
		\item The function $q(x)$ in the Helmholtz equation is taken to be constant zero.
		\item The ground truth source function $u_s$ is defined as
		\myen{
			u_{s}(x) = 0.5\exp(-300(x-0.4)^{2}) + 0.5\exp(-300(x-0.6)^{2}).
		}
	\end{itemize}
	
	According to the studies presented in \cite{Lim2016Thesis}, for this simple one-dimensional case,
	we will not take a recursive strategy but combine instead all data together with
	the forward operator denoted by $H$ and defined by
	$H = ( H_{\kappa_{1}}, H_{\kappa_{2}}, \cdots, H_{\kappa_{100}} )^{T}$.
	Based on these settings, we provide some basic theoretical properties of the prior and posterior sampling functions as follows:
	\begin{itemize}
		\item The prior probability measure for $u_s$ is Gaussian with the covariance operator $\mathcal{C}_{0}^{K}(\lambda)$
		with $\lambda\in\mathbb{R}^{+}$. According to Theorem 12 illustrated in \cite{Dashti2014}, we know that if $u_s$ is drawn from the
		prior measure, and then the following holds
		\myen{
			u_s \in W^{t,2}(\Omega) \quad \text{for }t < \frac{1}{2}, \quad \text{and} \quad
			u_s \in C^{0,t}(\Omega) \quad \text{for }t<\frac{1}{2},
		}
		where $W^{t,2}(\Omega)$ is the usual Sobolev space with $t$ times derivative belonging to $L^{2}(\Omega)$, and
		$C^{0,t}$ is the conventional H\"{o}lder space.
		\item For Algorithm \ref{algComplexVBG}, every posterior mean estimate $u_k$ has the following form:
		\myen{
			u_{k} = (\tau_{k}H^{*}H + \mathcal{C}_{0}(\lambda_{k})^{-1})^{-1}\tau_{k}H^{*}d.
		}
		Given that $H$ maps a function in $L^{2}(\Omega)$ to $\mathbb{R}^{200}$, we know that $H^{*}d$ is at least a
		function belonging to $L^{2}$. Considering the specific choices of $\mathcal{C}_{0}$,
		we have $u_k \in W^{2, 2}(\Omega)$. For Algorithm \ref{algComplexVBL}, we can derive similar conclusions.
	\end{itemize}
	
	\begin{remark}
		By employing the ``Bayesianize-then-discretize'' method, we can analyze the prior and posterior sampling functions rigorously.
		It is one of the advantages of employing our proposed infinite-dimensional VBM.
	\end{remark}
	
	\textbf{Gaussian noise case:}
	Let $d^{\dag}$ be the data without noise. Then, we construct noisy data by setting $d  = d^{\dag} + \sigma \xi$
	with $\sigma = 10^{-3}$ and $\xi$ is a random variable sampled from the standard normal distribution.
	
	
	\begin{figure}[htbp]
		\centering
		\includegraphics[width=0.8\textwidth]{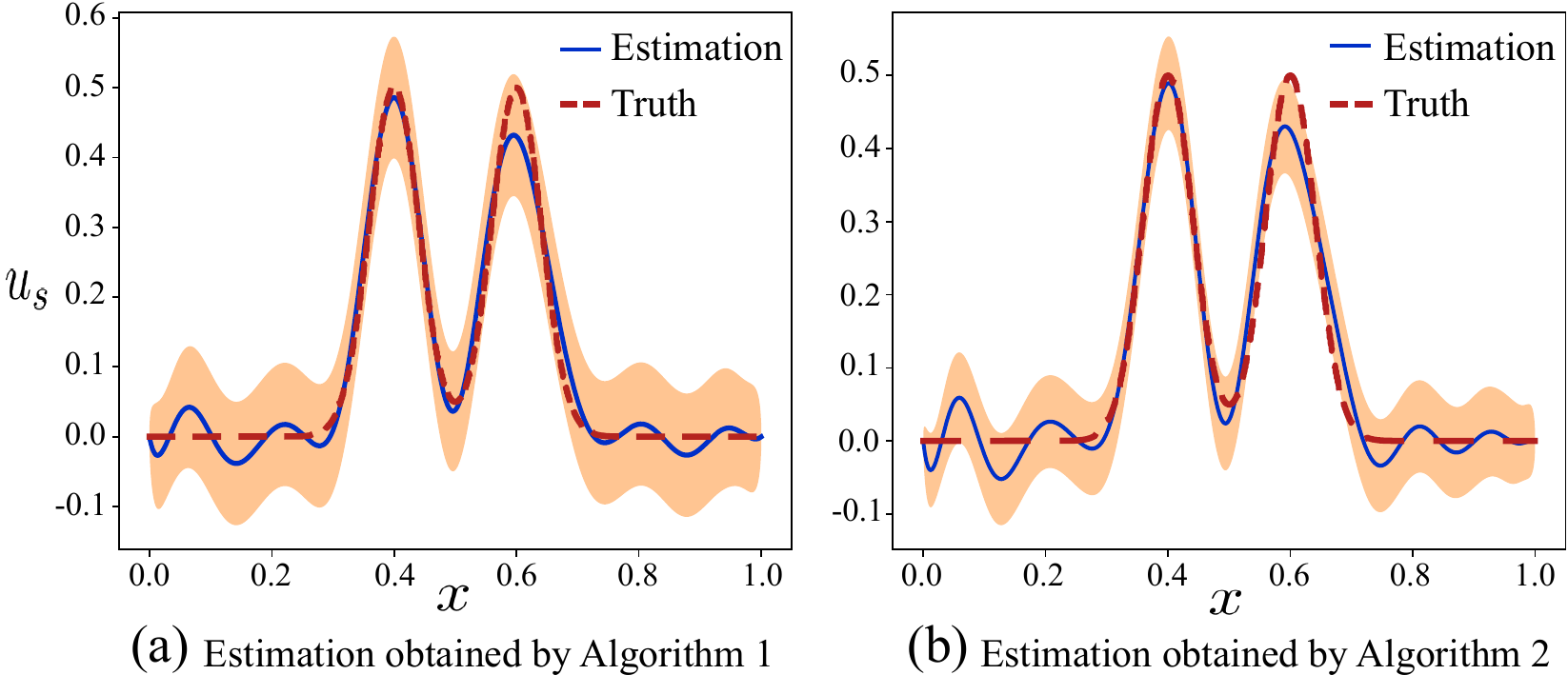}\\
		\vspace*{-0.6em}
		\caption{The truth and estimated functions when the data are polluted by Gaussian noise.
			(a): the estimated function obtained by Algorithm \ref{algComplexVBG} is denoted by the blue solid line,
			and the truth is denoted by the red dashed line;
			(b): the estimated function obtained by Algorithm \ref{algComplexVBL} is denoted by the blue solid line,
			and the truth is denoted by the red dashed line.
			In both plots the shaded areas represent the pointwise mean plus and minus two standard
			deviations from the mean (corresponding roughly to the $95\%$ confidence region).}\label{Ex1confidG}
	\end{figure}
	
	\vskip -0.2 cm
	\begin{figure}[htbp]
		\centering
		\includegraphics[width=0.8\textwidth]{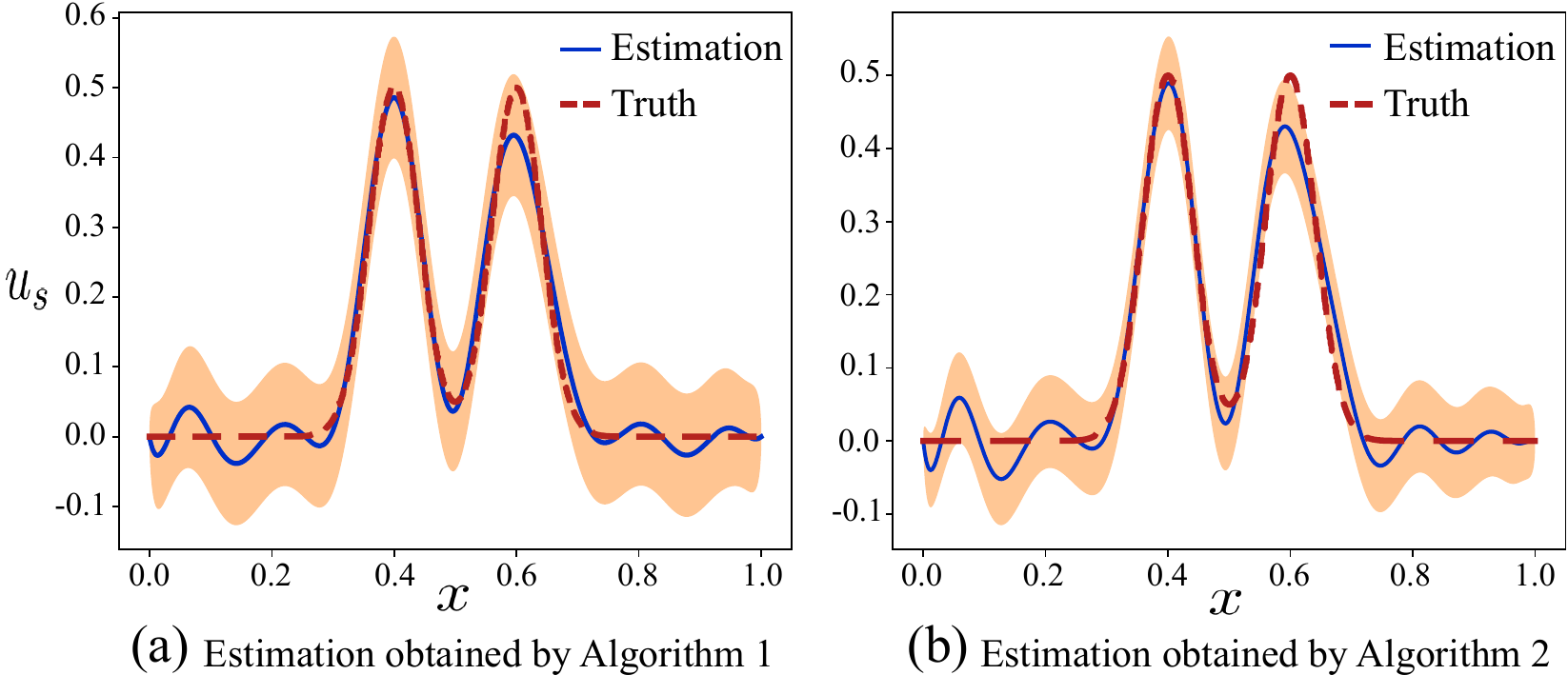}\\
		\vspace*{-0.6em}
		\caption{Relative errors of the estimated means in the $L^{\infty}$-norm of Algorithms \ref{algComplexVBG} and \ref{algComplexVBL}.}\label{Ex1DataL}
	\end{figure}
	
	In Figure \ref{Ex1confidG}, we depict the truth and estimated sources obtained by Algorithms \ref{algComplexVBG}
	and \ref{algComplexVBL}, respectively. Visually, both algorithms provide reasonable results.
	In addition, we demarcate the $95\%$ confidence region by the shaded area to display the uncertainties estimated by these two algorithms.
	The truth falls entirely into the confidence region given by Algorithm \ref{algComplexVBG},
	and the truth lies mostly within the confidence region given by Algorithm \ref{algComplexVBL}.
	This may indicate that for the Gaussian noise case, Algorithm \ref{algComplexVBG} can provide a more reliable estimation,
	which is in accordance with our assumptions.
	
	To give a more elaborate comparison, we present the relative errors of the estimated means in the $L^{\infty}$-norm of
	the two algorithms in Figure \ref{Ex1DataL}.
	The relative error of the conditional mean estimate used here is defined as follows
	\begin{align*}
	\text{relative error} = \|u - u_s\|_{L^{\infty}}/\|u_s\|_{L^{\infty}},
	\end{align*}
	where $u$ is the estimated function generated by our algorithm and $u_s$ is the true source function.
	The blue solid line and orange dashed line denote the relative errors obtained by Algorithms \ref{algComplexVBG} and \ref{algComplexVBL},
	respectively. Obviously, these two algorithms can provide comparable results after convergence.
	However, Algorithm \ref{algComplexVBG} converges much faster than Algorithm \ref{algComplexVBL},
	which is reasonable because the weight parameters used for detecting impulsive noises may reduce the convergence speed.
	
	The parameter $\tau$ given by Algorithm \ref{algComplexVBG} provides an estimate of the noise variance
	through $\sigma = \sqrt{\tau^{-1}}$. The true value of $\sigma$ is $0.001$ in our numerical example.
	To generate a repeatable results, we specify the random seeds in numpy to some certain numbers by
	$\text{numpy.random.seed(i)}$ with $i$ specified as some designated integers.
	The estimated $\sigma$ is equal to $0.000953, 0.001101, 0.001022, 0.001003$, and $0.001041$
	when the random seeds are specified as $1,2,3,4$, and $5$, respectively, thereby illustrating the
	effectiveness of our proposed algorithm.
	
\begin{figure}[htbp]
		\centering
		\includegraphics[width=0.4\textwidth]{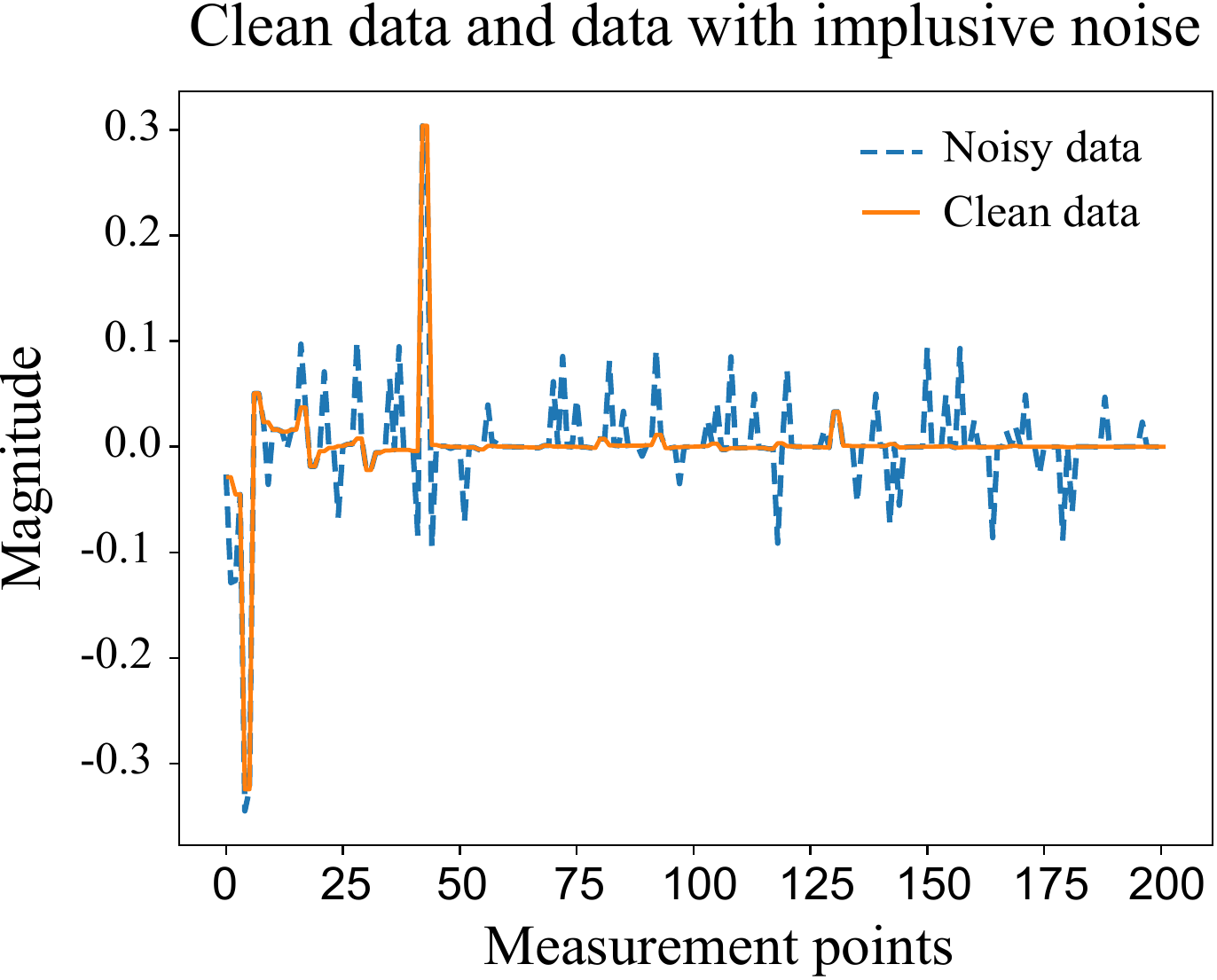}\\
		\vspace*{-0.6em}
		\caption{Clean and noisy data. The orange solid line represents the clean data, and
			the blue dashed line represents the data with impulsive noise.}\label{Ex1NoiseL}
	\end{figure}

	
	\textbf{Laplace noise case:} As for the Gaussian noise case, let $d^{\dag}$ be the noise-free measurement.
	The noisy data are generated as follows:
	\myen{
		d_{i} = \left\{\begin{aligned}
			& d_{i}^{\dag}, \qquad\quad \text{with probability }1-r, \\
			& d_{i}^{\dag} + \epsilon\xi, \quad \text{with probability }r,
		\end{aligned}\right.
	}
	where $\xi$ follows the uniform distribution $U[-1,1]$, and $(\epsilon, r)$ controls the noise pattern,
	$r$ is the corruption percentage, and $\epsilon$ is the corruption magnitude.
	In the following, we take $r = 0.5$ and $\epsilon = 0.1$.
	We plot the clean and noisy data in Figure \ref{Ex1NoiseL},
	which illustrates that the clean data are heavily polluted.
	
	\vskip -0.2 cm
	\begin{figure}[htbp]
		\centering
		\includegraphics[width=0.8\textwidth]{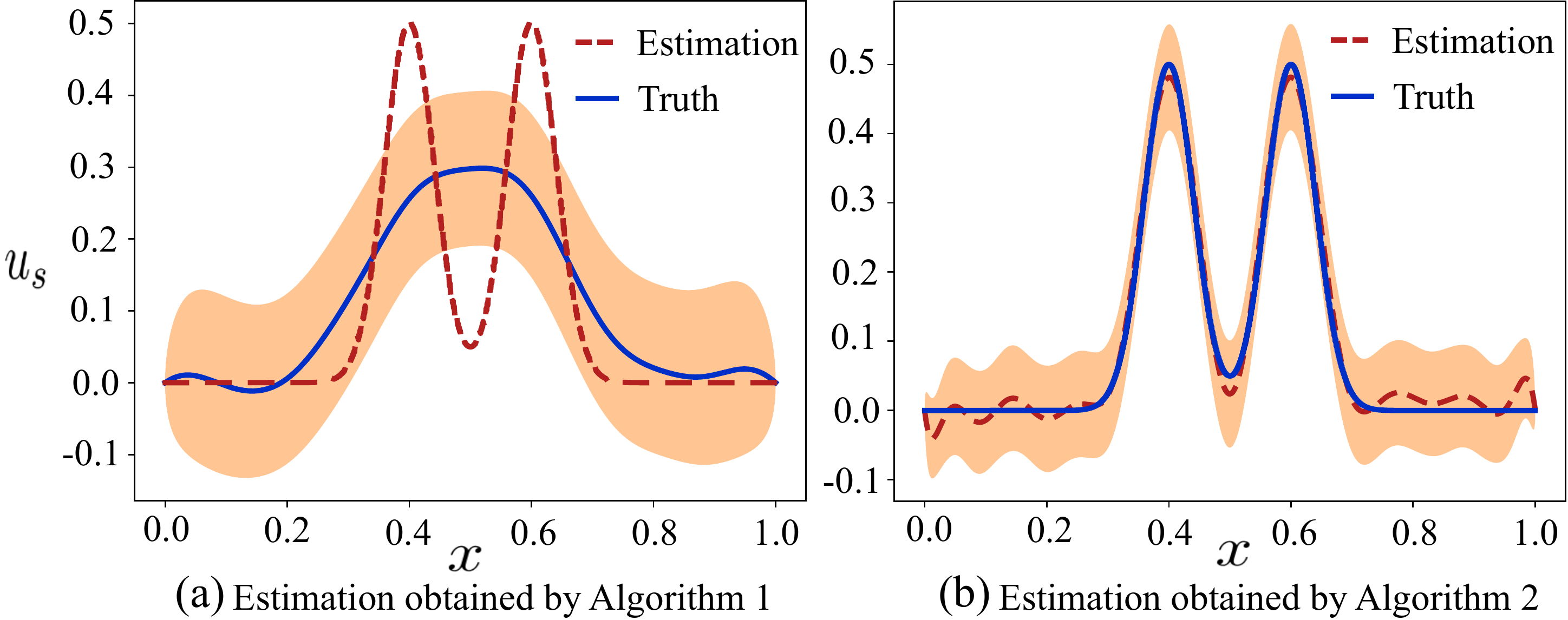}\\
		\vspace*{-0.6em}
		\caption{The truth and estimated functions when the data are polluted by impulsive noise.
			(a): The estimated function obtained by Algorithm \ref{algComplexVBG} is denoted by the blue solid line,
			and the truth is denoted by the red dashed line;
			(b): The estimated function obtained by Algorithm \ref{algComplexVBL} is denoted by the blue solid line,
			and the truth is denoted by the red dashed line.
			The shaded areas in both panels represent the pointwise mean plus and minus two standard
			deviations from the mean (corresponding roughly to the $95\%$ confidence region).}\label{Ex1confidL}
	\end{figure}
	
	\begin{figure}[htbp]
		\centering
		\includegraphics[width=0.9\textwidth]{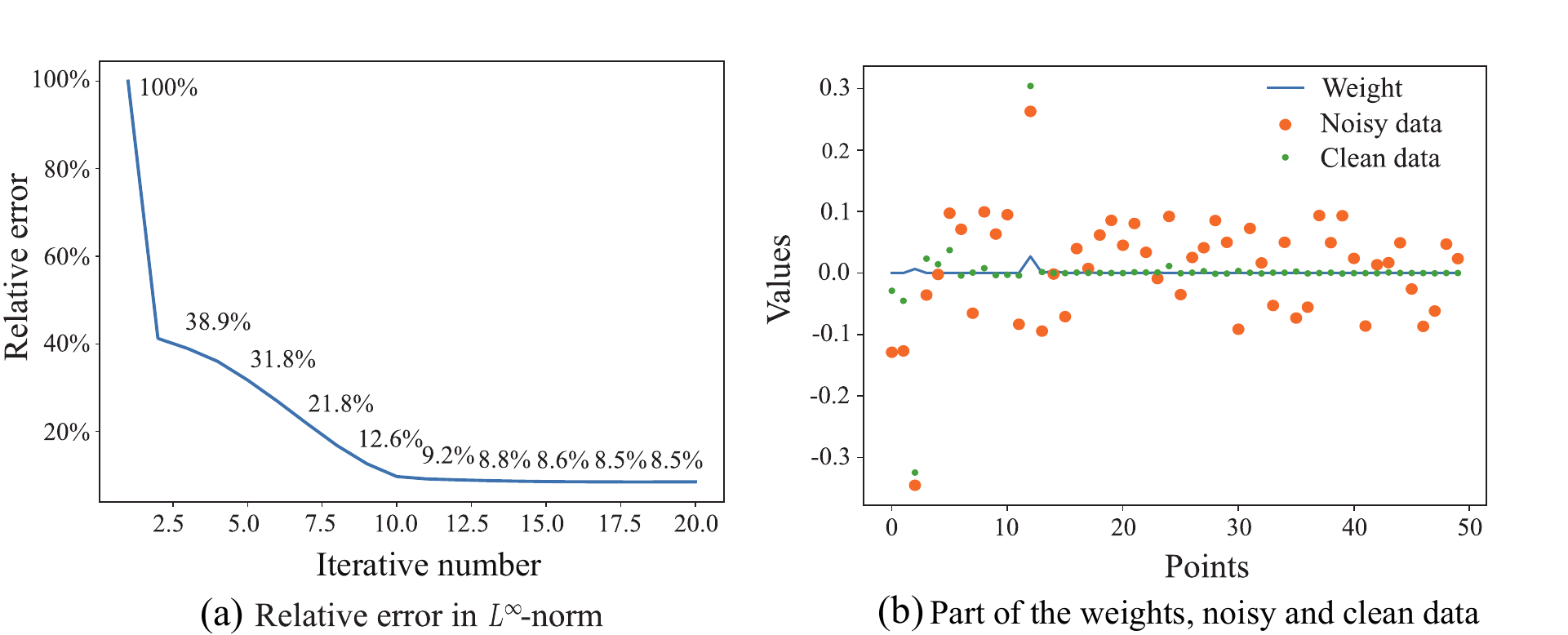}\\
		\vspace*{-0.6em}
		\caption{(a): Relative errors in the $L^{\infty}$-norm obtained by Algorithm 2;
			(b): Weight, noisy, and clean data at the data points with impulsive noise
			(only points with impulsive noise, not all points).} \label{Ex1ErrorL}
	\end{figure}
	In Figure \ref{Ex1confidL}, we show the estimated functions obtained by Algorithms \ref{algComplexVBG}
	and \ref{algComplexVBL} in the left and right panels, respectively.
	Obviously, based on the Gaussian noise assumption, Algorithm \ref{algComplexVBG} cannot provide a reasonable estimate, and
	the estimated confidence region may be unreliable. However, based on the Laplace noise assumption, Algorithm \ref{algComplexVBL}
	provides an accurate estimate. Given that Algorithm \ref{algComplexVBG} fails to converge to a reasonable estimation, we only
	provide the relative errors in the $L^{\infty}$-norm of Algorithm \ref{algComplexVBL} on the left panel of Figure \ref{Ex1ErrorL}.
	From these relative errors, we can find that Algorithm \ref{algComplexVBL} rapidly converges even if
	the data are heavily polluted by noise.
	The right panel of Figure \ref{Ex1ErrorL} plots the noisy and clean data points at those data points where noises are added.
	We plot the weight vector at the corresponding data points.
	From this figure, we can clearly see that the elements of the weight vector are all with small values,
	which is in accordance with our theory. The weight vectors at the noisy data points are adjusted to small values
	during the iteration. This reveals the outlier removal mechanism of Algorithm \ref{algComplexVBL}.

	\subsection{Two-dimensional ISP}\label{TwoDISPSection}
	In this subsection, we solve the two-dimensional ISP.
	Directly computing the covariance operator for the two-dimensional problem is difficult due to
	the large memory requirements and computational inefficiency.
	Here, we employ a simple method that employs a rough mesh approximation to compute the covariance.
	The source function $u_s$ can be expanded under basis functions as follows:
	\mye{
		u_s(x) = \sum_{i=1}^{\infty}u_{si}\varphi_{i}(x).
	}
	Given that these basis functions can be taken as the finite element basis, the source function can be approximated as
	\mye{\label{Ex2ApproUs}
		u_s(x) \approx \sum_{i=1}^{N_t}u_{si}\varphi_{i}(x).
	}
	The covariances involved in Algorithms \ref{algComplexVBG} and \ref{algComplexVBL} are all computed
	by taking a small $N_t$ in (\ref{Ex2ApproUs}). For many applications such as medical imaging,
	we may compute the operator $H^{*}H$ (not depending on the source function) with a small $N_t$
	before the inversion. To evaluate accurately as the wavenumber increases, we compute the
	mean function by gradient descent with a fine mesh discrete PDE solver and then project the source function
	to the rough mesh for computing variables relying on the covariance operators.
	
	Unlike the one-dimensional case, we employ the sequential method used in \cite{Bao2015IP} 
	that provides a more stable recovery for multi-frequency inverse problems.
	Specifically speaking, for $0=\kappa_0 < \kappa_1 < \cdots < \kappa_{N_f} < \infty$ and each problem
	$d_{\kappa_i} = H_{\kappa_i}(u_s) + \epsilon_{\kappa_i}$($i=1,\cdots,N_f$), we assume the prior measure is
	$\mu_{0i}^{u,\lambda} = \mathcal{N}(\bar{u}_{i-1},\mathcal{C}_0^{K}(\lambda))$ with $\bar{u}_{i-1}$ denoting the conditional mean
	estimate when the wavenumber is $\kappa_{i-1}$ ($\bar{u}_0$ is assumed to be some initial guess $u_s^0$).
	We denote $\|\cdot\|_{\mathcal{C}_0^K(\lambda)}$ as the Cameron-Martin norm corresponding to
	the Gaussian measure $\mathcal{N}(0,\mathcal{C}_0^{K}(\lambda))$.
	For the Gaussian noise case with $i=1,2,\cdots, N_f$, we have the following Bayesian formula
	\begin{align}\label{seql1}
	\frac{d\mu^{i}}{d\mu_{0i}}(u, \lambda, \tau) \propto \exp\Big(-\frac{\tau}{2}\|H_{\kappa_i}(u) - d_{\kappa_i}\|^2\Big),
	\end{align}
	where $\mu_{0i}(du, d\lambda, d\tau) = \mu_{0i}^{u,\lambda}(du)\mu_0^\lambda(d\lambda)\mu_0^{\tau}(d\tau)$ with
	$\mu_0^\lambda, \mu_0^{\tau}$ are defined as in Subsection \ref{SectionLinearIPG} and
	$\mu^i$ is the posterior measure when wavenumber is equal to $\kappa_i$.
	The posterior measure $\mu^{\kappa_{N_f}}$ will be employed to quantify the uncertainties of the final estimate.
	For a similar sequential formulation as above, we refer to Subsection 6.4.1 in \cite{Lim2016Thesis}.
	It is not hard to formulate a sequential approach for the Laplace noise case. The details are omitted for content conciseness.
	The iteration details are presented in Algorithm \ref{algComplexConPro}.
	\begin{algorithm}
		\caption{VBM for two-dimensional ISP with multi-frequencies}
		\label{algComplexConPro}
		\begin{algorithmic}
			\STATE {1: Give an initial guess of the unknown source $u_s$, denoted by $u_s^0$.
			}
			\STATE {2: For $i$ from $1$ to $N_f$ (iterate from low wavenumber to high wavenumber)
			}
			\STATE {3: \qquad Specify the prior measure of $u_s$ as
				$\mu_{0i}^{u,\lambda} = \mathcal{N}(u_s^{i-1},\mathcal{C}_0^{K}(\lambda))$.
				Running itera- \\
				\qquad\quad\,\! tions of Algorithms \ref{algComplexVBG} or \ref{algComplexVBL} for $k$
				until some stopping criterion is satisfied. \\
				\qquad\quad\,\! For $k=1$, rough approximate of $H$ and source is employed; For $k>1$, the \\
				\qquad\quad\,\! gradient descent method is employed to solve
				\begin{align*}
				u_s^{k}:=\argmin_{u_s} \frac{\tau_k}{2}\|H_{\kappa_i}(u_s) - d_{\kappa_i}\|^2 +
				\|u_s- u_s^{i-1}\|_{\mathcal{C}_0^{K}(\lambda_k)}^2,
				\end{align*}
				\qquad\quad\,\! which generate a conditional mean estimate on a fine mesh. In all of the \\
				\qquad\quad\,\! iterations, rough approximate Hessian has been used to update distributions \\
				\qquad\quad\,\! of hyper-parameters $\lambda$,
				$\tau$ (Algorithm \ref{algComplexVBG}) or $w$ (Algorithm \ref{algComplexVBL}).
			}
			\STATE {5: End for
			}
			\STATE {6: Return the approximate probability measure $\nu$.
			}
		\end{algorithmic}
	\end{algorithm}
	In the following, when we say that Algorithm \ref{algComplexVBG} is employed, we actually means
	that Algorithm \ref{algComplexConPro} is employed in combination with Algorithm \ref{algComplexVBG}.
	Similarly, when we say that Algorithm \ref{algComplexVBL} is employed, we mean that
	Algorithm \ref{algComplexConPro} is employed in combination with Algorithm \ref{algComplexVBL}.
	
	\begin{remark}\label{simplemethodiss}
		It should be pointed out that the simple ``rough mesh approximation'' method employed in Algorithm \ref{algComplexConPro}
		is only applicable to problems with a simple form
		(e.g., a localized source) on simple geometry. And this method can be hardly employed for dealing with more complex
		problems in three-dimension or even in two-dimension where a large $N_t$ is needed
		(e.g., high-resolution recovery with data of high wavenumbers).
		Our aim is to give an illustration of the proposed method. 	
		For more advanced techniques designed for large-scale problems, \cite{Tan2013IP} can be referred to, which
		provides a scalable discretize method for the infinite-dimensional Bayesian approach with linear approximations.
		The fully nonlinear case has been investigated by using a stochastic Newton MCMC method in \cite{Petra2014SIAM}.
		Then, Metropolize-then-discretize and discretize-then-Metropolize have been analyzed carefully for large-scale problems \cite{Thanh2016IPI}.
		In 2019, an approximate sampling method based on some randomized MAP estimates has been investigated in detail \cite{Wang2019SISC}.
		All these studies provide valuable ideas for designing algorithms of large-scale inverse problems.
		For more studies in this direction, we refer to \cite{Bui2014IP,Isaac2015JCP,Martin2012SISC}.
	\end{remark}
	
	
	\begin{remark}\label{Ex2explain}
		In Algorithm \ref{algComplexConPro}, we use approximations on a rough mesh for the first iteration of every wavenumber,
		which may provide an initial inaccurate adjustment for the parameters employed in
		Algorithms \ref{algComplexVBG} and \ref{algComplexVBL}.
		In our numerical experiments, we only take three iterations for the third step to obtain an estimation.
	\end{remark}
	
	
	\begin{remark}\label{ConExMCMC}
		To employ sampling-type methods such as the MCMC algorithm, researchers often parameterize the unknown source function
		carefully to reduce the dimension, e.g., assume that the sources are point sources, then parameterize the source function
		by numbers, locations, and amplitudes \cite{Engel2019IP}. For employing MCMC algorithm \cite{Cotter2013SS,Feng2018SIAM} in our setting,
		the computational complexity is unacceptable for two reasons:
		Calculation with many wavenumbers are needed for multi-frequency problems and
		a large number of samples need to be generated for each wavenumber; For each problem (\ref{seql1}),
		we did not assume any parametric form of the source function which makes the parameters of source equal to the dimension of
		the discretization (much more parameters than the usually used parametric form).
		However, the proposed Algorithm \ref{algComplexConPro} only takes several times of computational time
		compared with the classical iterative algorithms \cite{Bao2015IP,Bao2015SIAM,Isakov2018SIAM}
		to provide estimations of uncertainties.
	\end{remark}
	
	Before going further, we list the specific choices for some parameters introduced in Section \ref{SectionTypical} as follows:
	\begin{itemize}
		\item The operator $\mathcal{C}_{0}$ is chosen as $(-\Delta + \text{Id})^{-2}$. Here, the Laplace operator is
		defined on $\Omega$ with the zero Dirichlet boundary condition.
		\item Take the discrete truncate level $N_t = 1681$ and the number of measurement points $N_m = 200$.
		The basis functions $\{\varphi_{j}\}_{j=1}^{\infty}$ are specified as second-order
		finite element basis functions.
		\item For Algorithm \ref{algComplexConPro} combined with Algorithm \ref{algComplexVBG},
		the wavenumber series are specified as
		$\kappa_{j} = j \text{ with } j=1,3,5,\cdots,35.$
		For Algorithm \ref{algComplexConPro} combined with Algorithm \ref{algComplexVBL},
		the wavenumber series are specified as
		$\kappa_{j} = j \text{ with } j=1,2,3,\cdots,35.$
		\item The scatterer function $q(x)$ is defined as follows:
		\myen{
			q(x_1,x_2) = & 0.3(4-3x_1)^2e^{\left( -9(x_1-1)^2 - 9(x_2-2/3)^2 \right)} \\
			& - \left( 0.6(x_1-1) - 9(x_1-1)^3 - 3^5(x_2-1)^5 \right)e^{\left( -9(x_1-1)^2 - 9(x_2-1)^2 \right)} \\
			& - 0.03 e^{-9(x_1-2/3)^2 - 9(x_2-1)^2},
		}
		which is the function used in Subsection 2.6 in \cite{Bao2015IP}.
		\item The true source function $u_s$ is defined as follows:
		\myen{
			u_s(x) = 0.5e^{-100((x_1-0.7)^2 + (x_2-1)^2)} + 0.3e^{-100((x_1-1.3)^2+(x_2-1)^2)}.
		}
		\item To avoid the inverse crime, a mesh with mesh number $125000$ is employed for generating the data.
		For the inversion, two types of meshes are employed: a mesh with mesh number $28800$
		is employed when the wavenumbers are below $20$, and a mesh with mesh number $41472$
		is employed when the wavenumbers are greater than $20$.
	\end{itemize}
	
	
	\textbf{The case of Gaussian noise: }
	Let $d^{\dag}$ be the data without noise.
	The synthetic noisy data $d$ are generated by
	$d_{j} = d^{\dag} + \sigma \xi_{j}$,
	where $\sigma = \max_{1\leq j\leq N_m}\{|d^{\dag}_{j}|\}L_{\text{noise}}$ with
	$L_{\text{noise}}$ denoting the relative noise level and $\xi_{j}$ denoting the standard normal random variables.
	In our experiments, we take $L_{\text{noise}} = 0.05$, that is $5\%$ of noises are added.
	\begin{figure}[htbp]
		\centering
		\includegraphics[width=0.75\textwidth]{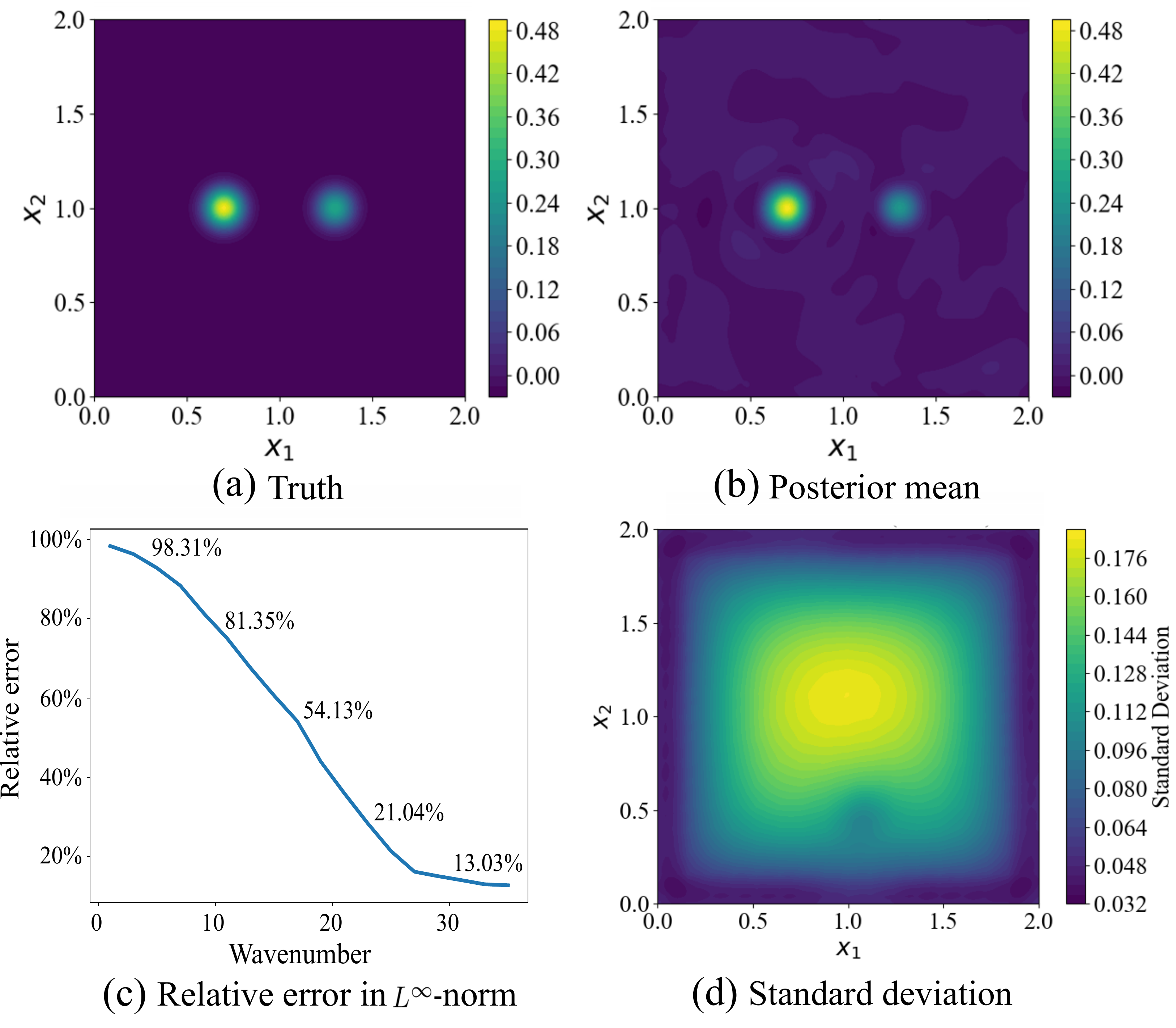}\\
		\vspace*{-0.6em}
		\caption{(a): The true source function;
			(b): The posterior mean estimate obtained by Algorithm \ref{algComplexVBG};
			(c): Relative error of the estimated means in $L^{\infty}$-norm obtained by Algorithm \ref{algComplexVBG};
			(d): Estimated standard deviation obtained by Algorithm \ref{algComplexVBG}.
		}\label{Ex2GaussianTrueEstimate}
	\end{figure}
	
	In Figure \ref{Ex2GaussianTrueEstimate}, we show the inference results obtained by Algorithm \ref{algComplexVBG}.
	We show the true source function on the top left and the posterior mean estimate on the top right.
	Visually, the estimate is similar to the truth, and only some small fluctuations in the background are observed.
	In the bottom left, we show the relative errors of the estimated means obtained by Algorithm \ref{algComplexVBG} as the wavenumber increases,
	which is in accordance with the results obtained by classical iterative approaches.
	In the bottom right, we show the estimated standard deviation obtained by Algorithm \ref{algComplexVBG}
	that quantifies the uncertainties of the posterior mean estimation.
	We see that the uncertainties are small on the boundary where data are collected.
	The areas with the largest uncertainties are in the middle, which is a reasonable result since
	that area can be recovered only when data generated by high wavenumbers are employed.

	\textbf{The case of Laplace noise: }
	For the Laplace noise case, let $d^{\dag}$ be the noise-free measurement.
	The noisy data are generated as
	\myen{
		d_{i} = \left\{\begin{aligned}
			& d_{i}^{\dag}, \qquad\quad \text{with probability }1-r, \\
			& d_{i}^{\dag} + \epsilon\xi, \quad \text{with probability }r,
		\end{aligned}\right.
	}
	where $\xi$ follows the uniform distribution $U[-1,1]$, $(\epsilon, r)$ controls the noise pattern,
	$r$ is the corruption percentage, and $\epsilon$ is the corruption magnitude defined by
	$\epsilon = \max_{1\leq j\leq N_m}\{|d^{\dag}_{j}|\}L_{\text{noise}}$ with
	$L_{\text{noise}}$ denoting the relative noise level.
	In our experiments, we take $L_{\text{noise}} = 1$ and $r=0.2 \text{ or } 0.5$.
	
	\begin{figure}[htbp]
		\centering
		\includegraphics[width=0.65\textwidth]{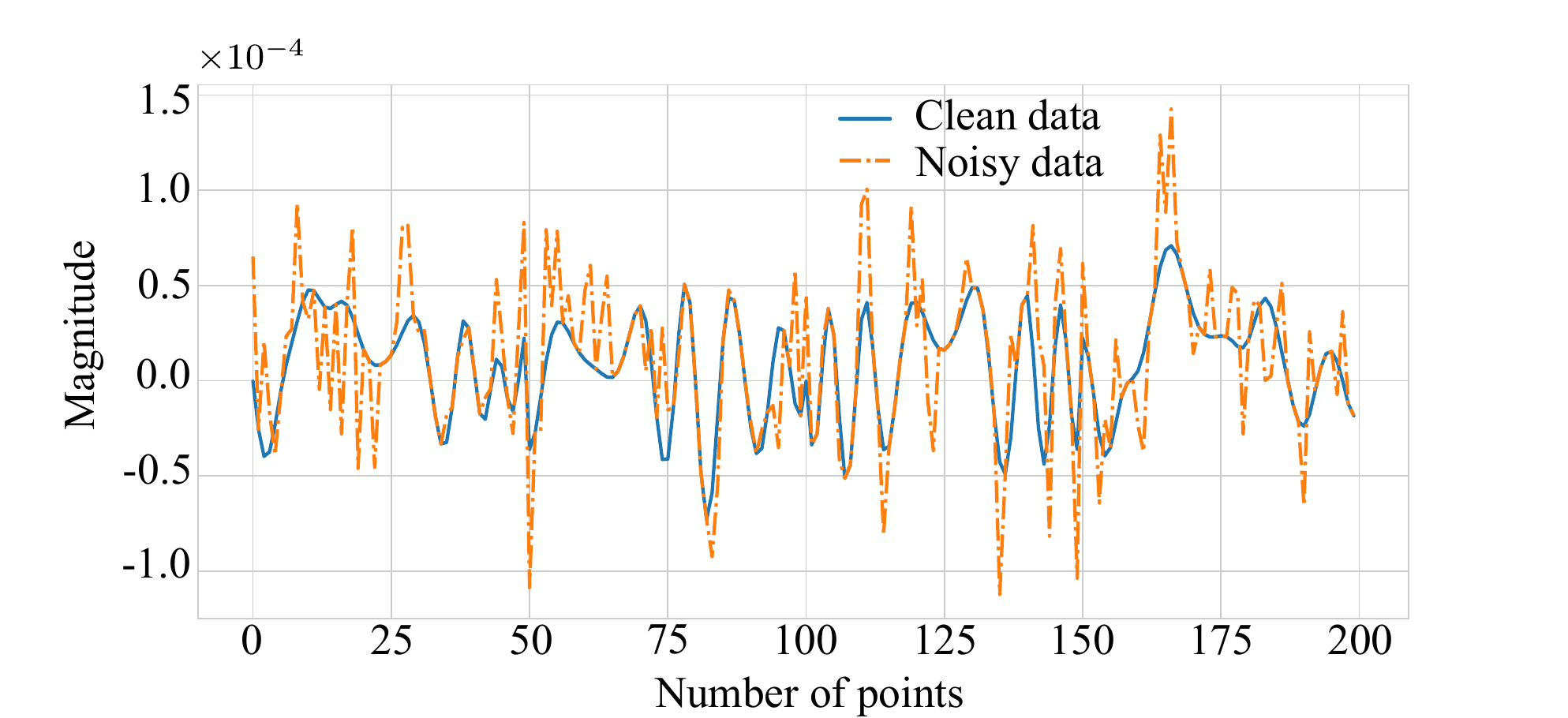}\\
		\vspace*{-0.6em}
		\caption{Clean and noisy data obtained when the wavenumber is $34$.
			The blue solid line represents the clean data, and
			the dashed orange line represents the noisy data with $r=0.5$.}\label{Ex2LaplaceNoisyData}
	\end{figure}
	\begin{figure}[htbp]
		\centering
		\includegraphics[width=0.95\textwidth]{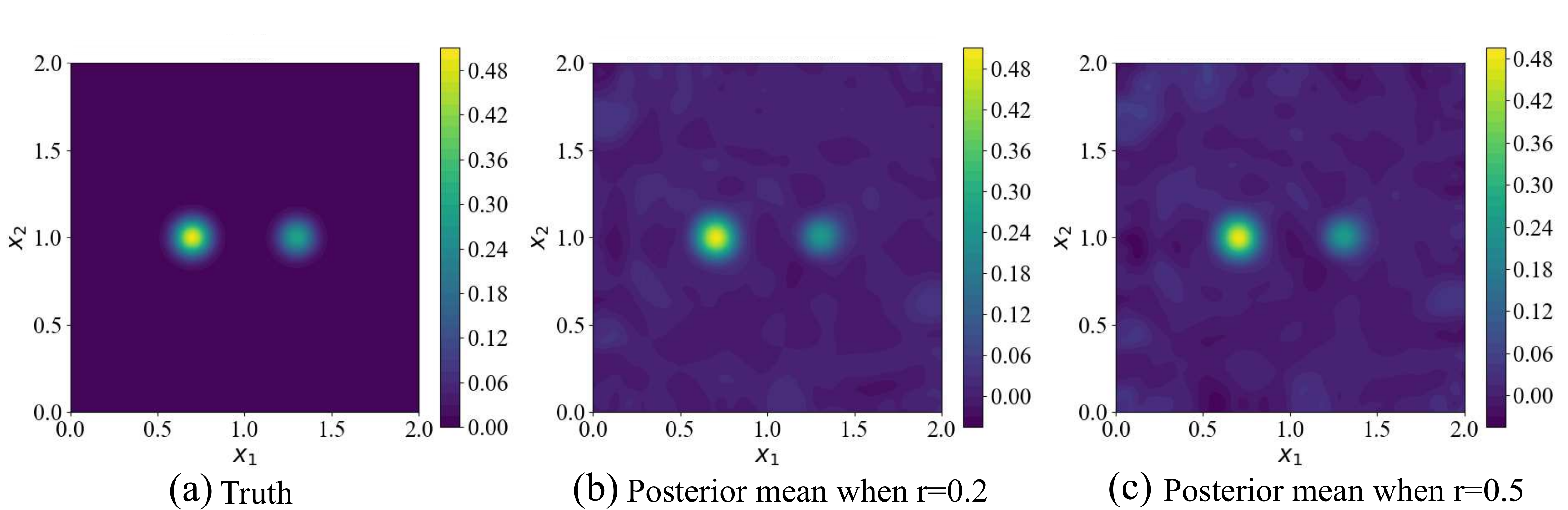}\\
		\vspace*{-0.6em}
		\caption{(a): The true source function;
			(b): The posterior mean estimate provided by Algorithm \ref{algComplexVBL}
			from noisy data with $r=0.2$ ($20\%$ of data are polluted);
			(c): The posterior mean estimate provided by Algorithm \ref{algComplexVBL}
			from noisy data with $r=0.5$ ($50\%$ of data are polluted).}\label{Ex2LaplaceTrueEstimate}
	\end{figure}
	\begin{figure}[htbp]
		\centering
		\includegraphics[width=0.8\textwidth]{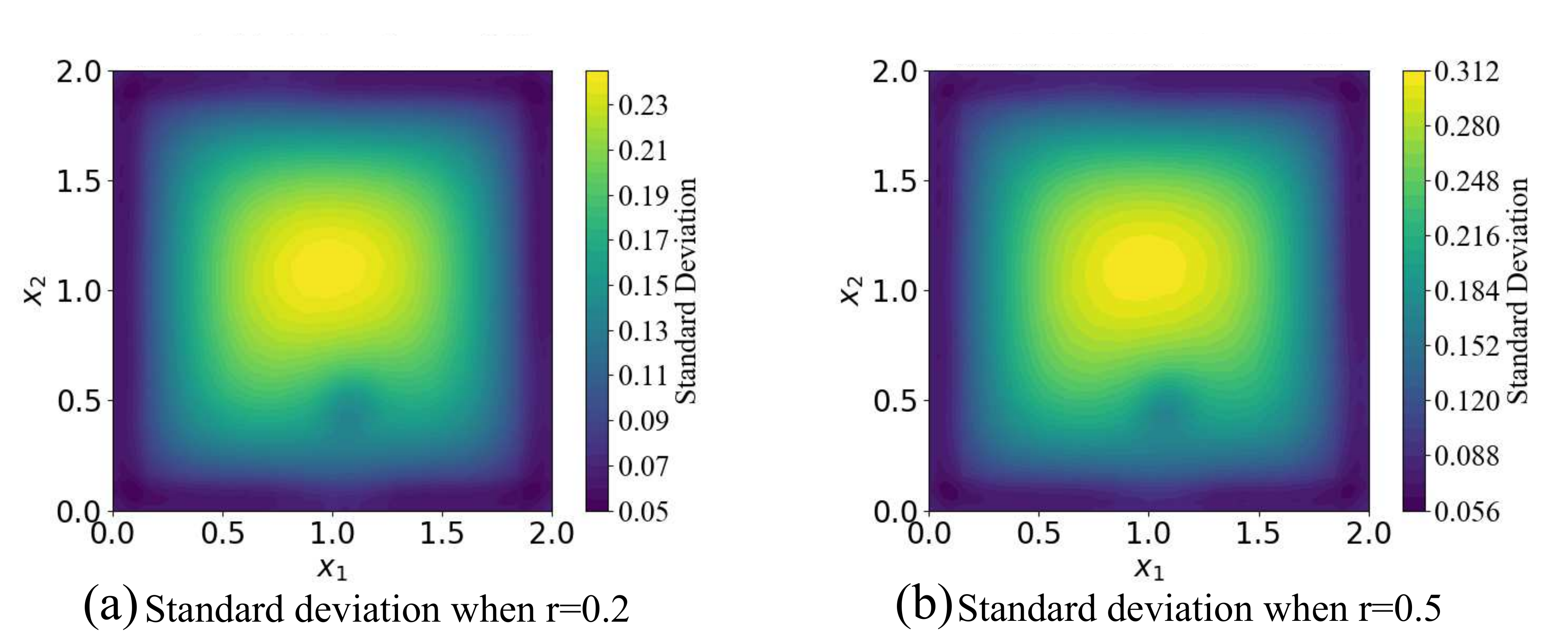}\\
		\vspace*{-0.6em}
		\caption{Standard deviation of the numerical solution obtained by Algorithm \ref{algComplexConPro}
			combined with Algorithm \ref{algComplexVBL}.
			(a): Estimated standard deviation when $r=0.2$ ($20\%$ of data are polluted);
			(b): Estimated standard deviation when $r=0.5$ ($50\%$ of data are polluted).}\label{Ex2LapCov2}
	\end{figure}
	\begin{figure}[htbp]
		\centering
		\includegraphics[width=0.75\textwidth]{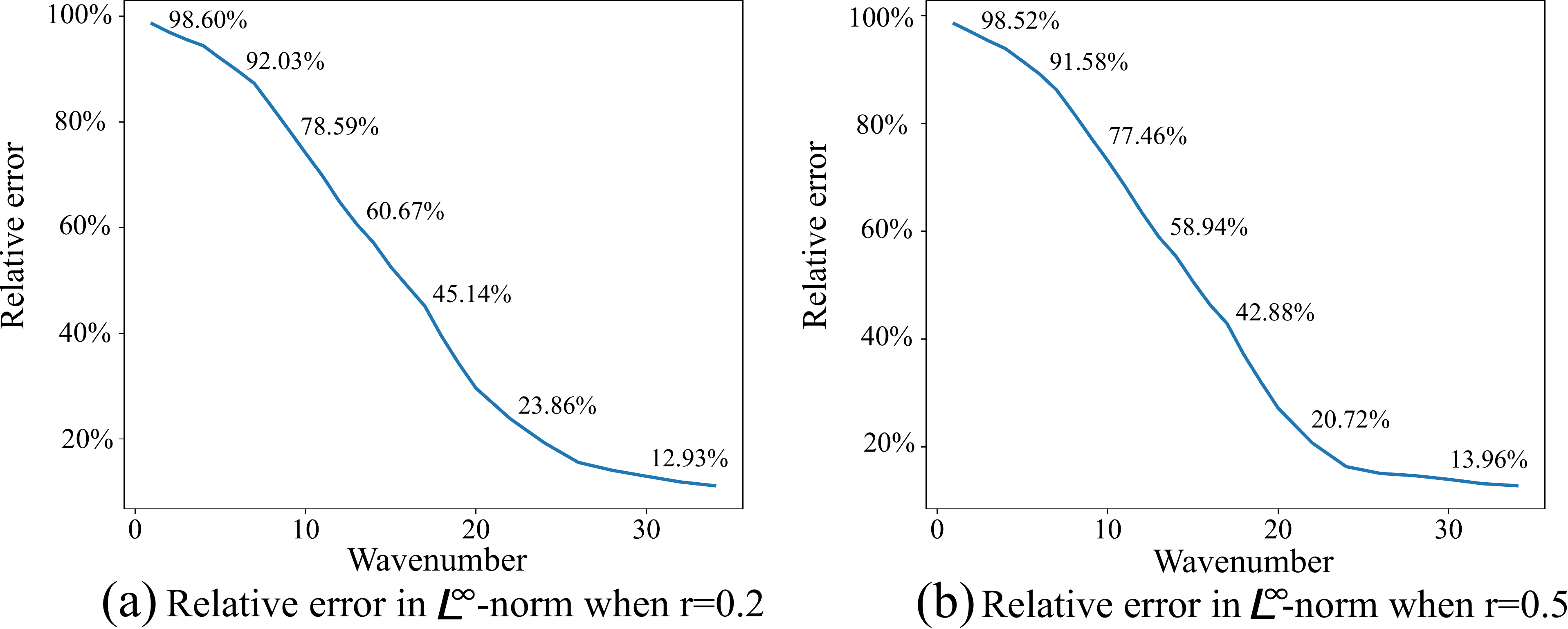}\\
		\vspace*{-0.6em}
		\caption{Relative errors of the estimated means in $L^{\infty}$-norm of Algorithm \ref{algComplexConPro} combined with Algorithm \ref{algComplexVBL}.
			(a): Relative errors for $r=0.2$; (b): Relative errors for $r=0.5$. }\label{Ex2LapRelativeError}
	\end{figure}
	
	
	The noisy and clean data when the wavenumber is $34$ and $r=0.5$ are shown in Figure \ref{Ex2LaplaceNoisyData}.
	Obviously, the data are heavily contaminated by noise.
	Figure \ref{Ex2LaplaceTrueEstimate} shows the true source function and the posterior mean estimates generated by Algorithm \ref{algComplexVBL}
	when $r=0.2$ and $r=0.5$ on the left, middle, and right panels, respectively.
	No essential differences can be observed between the posterior mean estimates when $r=0.2$ and $r=0.5$.
	However, the Bayes' method not only provides point estimates (e.g., posterior mean estimates) but also
	delivers the reliability of the  obtained estimations.
	Figure \ref{Ex2LapCov2} shows the standard deviations provided by Algorithm \ref{algComplexVBL} when $r=0.2$
	and $r=0.5$ on the left and right panels, respectively. The standard deviations are smaller when $r=0.2$, which is
	reasonable given that $80\%$ of the data are clean and only $50\%$ of the data are clean when $r=0.5$.
	Figure \ref{Ex2LapRelativeError} shows the relative errors in $L^{\infty}$-norm obtained by Algorithm \ref{algComplexVBL}
	with $r=0.2, 0.5$ on the left and right panels, respectively. Under both settings, the relative errors of the posterior mean
	estimates rapidly decrease.
	
	
	\begin{remark}
		The wavenumber series in the present paper are not chosen carefully in an optimal way.
		There are some studies focused on the strategies for selecting appropriate wavenumbers
		to give an accurate estimate under the framework of regularization methods
		for geophysical inverse problems \cite{Pan2019Geophysics}.
		Here, we choose more wavenumbers for the Laplace noise model based on a simple intuitive idea.
		More data are required when more hyper-parameters need to be inferred (The Laplace noise model has more parameters than
		the Gaussian noise model). 
	\end{remark}
	
	\section{Conclusion}
	In this paper, we have generalized the finite-dimensional mean-field approximate based variational Bayes' method (VBM)
	to infinite-dimensional space, which provides a mathematical foundation for applying VBM to the inverse problems of PDEs.
	A general theory for the existence of minimizers has been established, and by introducing the concept of reference
	probability measure, the mean-field approximate theory has been constructed for functions. The established general theory
	is then applied to abstract linear inverse problems with Gaussian and Laplace noise assumptions.
	Numerical examples for the inverse source problems of Helmholtz equations are investigated in details to highlight
	the effectiveness of the proposed theory and algorithms.
	
	There are numerous interesting problems that are worthy of being further investigated.
	Introducing a more reasonable setting of the intrinsic dimension will be meaningful.
	The recently published paper \cite{NIPS2019_9649} may provide some promising ideas. 
	For the infinite-dimensional Bayesian method with hyper-parameters, noncentered parameterization \cite{Agapiou2014SIAM} could be 
	a more appropriate choice. Using the proposed theory under the noncentered parameterization is 
	a problem worthy of being further investigated. 
	
	
	\section{Appendix}
	In this supplementary material, we provide all of the proof details for the lemmas and theorems presented in the main text. 
	
	\noindent\textbf{Proof of Lemma 2}
	\begin{proof}
		Let $\{\nu_n\}_{n=1}^{\infty} = \{ \prod_{i=1}^{M}\nu_{n}^{i} \}_{n=1}^{\infty}$ be a sequence of measures 
		in $\mathcal{C}$ that converges weakly to a probability measure $\nu_{*}$. 
		We want to show that $\nu_{*} \in \mathcal{C}$. Define 
		\mye{
			\nu_{*}^{i} := \int_{\prod_{j\neq i}\mathcal{H}_{j}}d\nu_{*}, \quad \text{for}\quad i=1,2,\cdots,M.
		} 
		Obviously, each $\nu_{*}^{i}$ is a probability measure. 
		Let $f_{i}$ be some bounded continuous function defined on $\mathcal{H}_{i}$ with $i=1,2,\cdots, M$. 
		Based on the definition of weak convergence, we obtain 
		\mye{\label{limit1}
			\int_{\prod_{j=1}^{M}\mathcal{H}_{j}}f_{i}d\nu_{n} \rightarrow \int_{\mathcal{H}_{i}}f_{i}d\nu_{*}^{i}, \quad
			\text{as} \,\,\, n\rightarrow\infty.
		}
		It should be noted that the left hand side of (\ref{limit1}) is equal to 
		\mye{
			\int_{\mathcal{H}_{i}}f_{i}d\nu_{n}^{i},
		} 
		and we find that each $\nu_{n}^{i}$ converges weakly to $\nu_{*}^{i}$. 
		Therefore, we find that $\nu_{*}^{i}$ belongs to $\mathcal{A}_{i}$. 
		Let $f$ be a bounded continuous function defined on $\prod_{j=1}^{M}\mathcal{H}_{j}$. 
		Then, it is a bounded continuous function for each variable. Based on the definition of weak convergence, we find that
		\mye{
			\int_{\prod_{j=1}^{M}\mathcal{H}_{j}}f d\nu_{n} \rightarrow \int_{\prod_{j=1}^{M}\mathcal{H}_{j}}f d\nu_{*},
		} 
		and
		\mye{
			\int_{\prod_{j=1}^{M}\mathcal{H}_{j}}f d\nu_{n} 
			& = \int_{\prod_{j=1}^{M}\mathcal{H}_{j}}f d\nu_{n}^{1}\cdots d\nu_{n}^{M}  
			\rightarrow \int_{\prod_{j=1}^{M}\mathcal{H}_{j}}f d\nu_{*}^{1}\cdots d\nu_{*}^{M},
		}
		when $n\rightarrow\infty$. Relying on the arbitrariness of $f$, we conclude that 
		$\nu_{*} = \prod_{j=1}^{M}\nu_{*}^{j}$, which completes the proof.  
	\end{proof}
	
	\noindent\textbf{Proof of Theorem 5}
	\begin{proof}
		From the proof of Lemma 2, we know that $\nu_{n}^{j}$ converges weakly to $\nu_{*}^{j}$
		for every $j = 1,2,\cdots,M$. According to $\nu_{n}^{j}\ll \nu_{*}^{j}$ for $j=1,2,\cdots,M$, we have
		\mye{\label{tvequ1}
			D_{\text{KL}}(\nu_{n}||\nu_{*}) & = \int\frac{d\nu_{n}}{d\nu_{*}}\log\bigg(\frac{d\nu_{n}}{d\nu_{*}}\bigg)d\nu_{*}
			= \sum_{j=1}^{M}\int\log\bigg( \frac{d\nu_{n}^{j}}{d\nu_{*}^{j}} \bigg)d\nu_{n}^{j}	\\
			& = \sum_{j=1}^{M}D_{\text{KL}}(\nu_{n}^{j} || \nu_{*}^{j}).
		}
		Using Lemma 2.4 proved in \cite{Pinski2015SIAMMA} and Lemma 22 shown in \cite{Dashti2014}, 
		we find that $\nu_{n}$ converges to $\nu_{*}$ in the total-variation norm.
		Combined with the above equality (\ref{tvequ1}), the proof is completed.	
	\end{proof}
	
	\noindent\textbf{Proof of Theorem 9}
	\begin{proof}
		For a fixed $j$, let $B\in\mathcal{M}(\mathcal{H}_{j})$, and $\nu_{n}^{j}\in\mathcal{A}_{j}$ be a sequence that converges 
		weakly to $\nu_{*}^{j}$ and  
		\mye{
			\frac{d\nu_{n}^{j}}{d\mu_{r}^{j}} = \frac{1}{Z_{nr}^{j}}\exp(-\Phi_{j}^{nr}(x_{j})).
		}
		Assuming that $\mu_{r}^{j}(B) = 0$ and by assumption (16) in the main text, we have 
		\myen{
			\nu_{n}^{j}(B) = \int_{B}\frac{1}{Z_{nr}^{j}}\exp(-\Phi_{j}^{nr}(x_{j}))\mu_{r}^{j}(dx_j) 
			= 0.
		}
		Define
		\mye{
			B_{m} = \{ x\in B \, | \, \text{dist}(x, B^{c}) \geq 1/m \},
		}
		and let $f_{m} > 0$ be a positive continuous function that satisfies
		\myen{
			f_{m}(x) = \left\{\begin{aligned}
				& 1, \quad x\in B_{m}, \\
				& 0, \quad x\in B^{c}.
			\end{aligned}\right.
		}
		Then, we have 
		\mye{
			\nu_{*}^{j}(B_{m}) \leq \int_{\mathcal{H}_{j}}f_{m}d\nu_{*}^{j} = 
			\lim_{n\rightarrow\infty}\int_{\mathcal{H}_{j}}f_{m}d\nu_{n}^{j} 
			\leq \lim_{n\rightarrow\infty} \nu_{n}^{j}(B) = 0,
		}
		and
		\mye{
			\nu_{*}^{j}(B) = \sup_{m}\nu_{*}^{j}(B_{m}) = 0,
		}
		based on the inner regular property of finite Borel measures. 
		Therefore, there exists a constant and a continuous function denoted by $Z_{r}^{j}$ and $\Phi_{j}^{r}(\cdot)$ such that 
		\mye{
			\frac{d\nu_{*}^{j}}{d\mu_{r}^{j}}(x_{j}) = \frac{1}{Z_{r}^{j}}\exp\big( -\Phi_{j}^{r}(x_{j}) \big).
		}
		To complete the proof, we should verify the almost surely positiveness of the right-hand side of the above equality.  
		Assume that $\frac{1}{Z_{r}^{j}}\exp\big( -\Phi_{j}^{r}(x_{j}) \big) = 0$ on a set $B\subset\mathcal{H}_{j}$ with 
		$\mu_{r}^{j}(B) > 0$. If $B\subset\mathcal{H}_{j}\backslash\sup_{N}T_{N}^{j}$, then it holds that $\mu_{r}^{j}(B) = 0$ by our 
		assumption. Therefore, $B\cap\sup_{N}T_{N}^{j}$ is not empty, and there exists a constant $\tilde{N}$ such 
		that for all $N \geq \tilde{N}$, $B\cap T_{N}^{j}$ is not empty. Denote $B_{N} = B\cap T_{N}^{j}$, 
		and then for a sufficiently large $N$, we have $\mu_{r}^{j}(B_{N}) \geq \frac{1}{2}\mu_{r}^{j}(B)$. 
		Let 
		\myen{
			B_{N}^{m} = \{ x\in B_{N} \, | \, \text{dist}(x, B_{N}^{c}) \geq 1/m \},
		}
		and define a function $g_{m}$ similar to $f_{m}$ with $B_{m}$ replaced by $B_{N}^{m}$. 
		Given that $\mu_{r}^{j}(B_{N}) = \sup_{m}\mu_{r}^{j}(B_{N}^{m})$, for a large enough $m$, we find that
		\myen{
			\mu_{r}^{j}(B_{N}^{m}) \geq \frac{1}{2}\mu_{r}^{j}(B_{N}) \geq \frac{1}{4}\mu_{r}^{j}(B) > 0.
		}
		By the definition of weak convergence, we have 
		\mye{
			\lim_{n\rightarrow \infty}\int_{\mathcal{H}_{j}}g_{m}(x) \frac{1}{Z_{nr}^{j}}\exp\big( -\Phi_{j}^{nr}(x) \big) d\mu_{r}^{j} = 
			\int_{\mathcal{H}_{j}}g_{m}(x)\frac{1}{Z_{r}^{j}}\exp\big( -\Phi_{j}^{r}(x) \big)d\mu_{r}^{j}.
		}
		The right hand side of the above equation is equal to $0$, but for a large enough $m$, the left hand side is positive and 
		the lower bound is 
		\mye{
			\frac{1}{4}\exp(-C_{N})\mu_{r}^{j}(B). 
		}
		This is a contradiction, and thereby the closedness of $\mathcal{A}_j$($j=1,\cdots,M$) have been proved.
		Combining the obtained results with the statements in Theorem 3, we obviously obtain the existence of a solution which 
		completes the proof.
	\end{proof}
	
	\noindent\textbf{Proof of Theorem 10}
	\begin{proof}
		Here, we focus on the deduction of formula (21) presented in the main text. 
		By inserting the prior probability measure into the Kullback-Leibler divergence between $\nu$ and $\mu$, 
		for each $i$ ($i=1,2,\cdots, M$) we find that 
		\myen{
			D_{\text{KL}}(\nu || \mu) & = \int_{\mathcal{H}}\log\bigg( \frac{d\nu}{d\mu_{r}} \bigg)
			- \log\bigg( \frac{d\mu_{0}}{d\mu_{r}} \bigg) - \log\bigg( \frac{d\mu}{d\mu_{0}} \bigg) d\nu \\ 
			& = \int_{\mathcal{H}} \bigg(-\sum_{j=1}^{M}\Phi_{j}^{r}(x_{j}) + \Phi^{0}(x) + \Phi(x) 
			\bigg) d\nu + \text{Const} \\
			& = \int_{\mathcal{H}_{i}} \bigg[ \int_{\prod_{j\neq i}\mathcal{H}_{j}} \bigg(\Phi^{0}(x) + \Phi(x)
			\bigg) \prod_{j\neq i}\nu^{j}(dx_{j}) \bigg] \nu^{i}(dx_{i})  \\
			& \quad
			- \int_{\mathcal{H}_{i}}\Phi_{i}^{r}(x_{i})\nu^{i}(dx_{i}) + \text{terms not related to }\Phi_{i}(x_{i}). 
		}
		For $i=1,2,\cdots,M$, let $\tilde{\nu}^{i}$ be a probability measure defined as follows:
		\mye{\label{defp1}
			\frac{d\tilde{\nu}^{i}}{d\mu_{r}^{i}} \propto \exp\bigg(-
			\int_{\prod_{j\neq i}\mathcal{H}_{j}} \bigg(\Phi^{0}(x) + \Phi(x)
			\bigg) \prod_{j\neq i}\nu^{j}(dx_{j}) \bigg).
		}
		By assumption (19) and (20) shown in the main text, we know that the right-hand side of (\ref{defp1}) is positive almost surely.
		Then, we easily know that the measures $\tilde{\nu}^{i}$ and $\mu_r^i$ are equivalent with each other.  
		Therefore, we obtain
		\mye{
			D_{\text{KL}}(\nu || \mu) & = -\int_{\mathcal{H}_{i}}\log\bigg( \frac{d\tilde{\nu}^{i}}{d\mu_{r}^{i}} \bigg)d\nu^{i}
			+ \int_{\mathcal{H}_{i}}\log\bigg( \frac{d\nu^{i}}{d\mu_{r}^{i}} \bigg)d\nu^{i} + \text{Const} \\
			& = D_{\text{KL}}(\nu^{i} || \tilde{\nu}^{i}) + \text{terms not related to $\nu^{i}$}.
		}
		Obviously, in order to attain the infimum of the Kullback-Leibler divergence, we should take
		$\nu^{i} = \tilde{\nu}^{i}$. Comparing formula (\ref{defp1}) with definition (14) in the main text, we notice that 
		the condition $\nu^{i} = \tilde{\nu}^{i}$ implies the following equality:
		\myen{
			\Phi_{i}^{r}(x_{i}) & = \int_{\prod_{j\neq i}\mathcal{H}_{j}}\bigg(\Phi^{0}(x) + 
			\Phi(x)\bigg) \prod_{j\neq i}\nu^{j}(dx_{j})
			+ \text{Const}, 
		}  
		which completes the proof.
	\end{proof}
	
	\noindent\textbf{Verify conditions in Theorem 10 for the linear inverse problem introduced in Subsection 3.1}\\
	\\
	At last, we provide a detailed verification of the conditions in Theorem 10 for the example employed in Subsection 3.1.
	As stated in Remark 14, we consider $\lambda' = \log\lambda$ and $\tau'=\log\tau$ as hyper-parameters. 
	For a sufficiently small $\epsilon>0$, taking $a_u(\epsilon, u):= \|u\|_{\mathcal{H}_u}^2$, 
	$a_{\lambda'}(\epsilon,\lambda') := \max\big\{-\lambda', \exp(\epsilon\exp(\lambda'))\big\}$ and 
	$a_{\tau'}(\epsilon,\tau') := \max\big\{-\tau', \exp(\epsilon\exp(\tau'))\big\}$, 
	then we try to verify conditions (19) and (20).
	In the following, the notation $C$ is a constant that may be different from line to line.
	In this example, we take $x_1 = u$, $x_2 = \lambda'$, and $x_3 = \tau'$. As shown in the main text, we have 
	\begin{align*}
	& \Phi^0(u,\lambda',\tau') = \frac{1}{2}\sum_{j=1}^{K}(u_j - u_{0j})^2(e^{\lambda'}-1)\alpha_j^{-1} - \frac{K}{2}\lambda', \\
	& \Phi(u, \lambda',\tau') = \frac{e^{\tau'}}{2}\|Hu-d\|^2 - \frac{N_d}{2}\tau'.
	\end{align*}
	With these preparations, we firstly verify 
	\begin{align}\label{cond1}
	T^{1}:=\sup_{u\in T_{N}^{u}}\sup_{\substack{\nu^{\lambda'}\in\mathcal{A}_{\lambda'}\\ \nu^{\tau'}\in\mathcal{A}_{\tau'}}}
	\int_{\mathbb{R}}\int_{\mathbb{R}}(\Phi^0 + \Phi)1_{A}(u,\lambda',\tau')\nu^{\lambda'}(d\lambda')\nu^{\tau'}(d\tau') < \infty.
	\end{align}
	Taking the specific expressions of $\Phi^0$ and $\Phi$ into (\ref{cond1}), we have 
	\begin{align}
	\begin{split}
	T^1 \leq C \sup_{u\in T_{N}^{u}}\sup_{\substack{\nu^{\lambda'}\in\mathcal{A}_{\lambda'}\\ \nu^{\tau'}\in\mathcal{A}_{\tau'}}} 
	\left( T^{11} + T^{12} + T^{13} + T^{14} \right),
	\end{split}
	\end{align}
	where 
	\begin{align*}
	& T^{11} = \int_{\mathbb{R}^{+}}\int_{\mathbb{R}}\frac{1}{2}\sum_{j=1}^{K}(u_j - u_{0j})^2 (e^{\lambda'}-1)\alpha_j^{-1}
	e^{-\Phi_{\tau'}^{r}(\tau')}e^{-\Phi_{\lambda'}^{r}(\lambda')}\mu_{r}^{\tau'}(d\tau')\mu_{r}^{\lambda'}(d\lambda'),  \\
	& T^{12} = \int_{\mathbb{R}^{-}}-\frac{K}{2}\lambda'e^{-\Phi_{\lambda'}^{r}(\lambda')}\mu_{r}^{\lambda'}(d\lambda'), \\
	& T^{13} = \int_{\mathbb{R}}\frac{e^{\tau'}}{2}\|Hu-d\|^2 e^{-\Phi_{\tau'}^{r}(\tau')}\mu_{r}^{\tau'}(d\tau'), \\
	& T^{14} = \int_{\mathbb{R}^{-}}-\frac{N_d}{2}\tau'e^{-\Phi_{\tau'}^{r}(\tau')}\mu_{r}^{\tau'}(d\tau'). 
	\end{align*}
	Because the techniques used for estimating these terms are similar, 
	we provide the estimates of $T^{13}$ as an example and omit the details for other terms. 
	Because $H$ is assumed to be a linear bounded operator, we have 
	\begin{align}
	\begin{split}
	T^{13} \leq & C \int_{\mathbb{R}} (e^{\epsilon e^{\tau'}} + 1)e^{-\Phi_{\tau'}^{r}(\tau')} \mu_{r}^{\tau'}(d\tau') \\
	\leq & C \int_{\mathbb{R}} \max(1, a_{\tau'}(\epsilon,\tau'))e^{-\Phi_{\tau'}^{r}(\tau')} \mu_{r}^{\tau'}(d\tau') < \infty. 
	\end{split}
	\end{align}
	Next, we need to estimate 
	\begin{align}\label{cond2}
	T^{2} := \sup_{\lambda'\in T_{N}^{\lambda'}}\sup_{\substack{\nu^{u}\in\mathcal{A}_{u}\\ \nu^{\tau'}\in\mathcal{A}_{\tau'}}}
	\int_{\mathcal{H}_u}\int_{\mathbb{R}}(\Phi^0 + \Phi)1_{A}(u,\lambda',\tau')\nu^{\tau'}(d\tau')\nu^{u}(du) < \infty.
	\end{align}
	Taking the specific expressions of $\Phi^0$ and $\Phi$ into (\ref{cond2}), we have
	\begin{align}
	T^2 \leq C \sup_{\lambda'\in T_{N}^{\lambda'}}\sup_{\substack{\nu^{u}\in\mathcal{A}_{u}\\ \nu^{\tau'}\in\mathcal{A}_{\tau'}}}
	\left( T^{21} + T^{22} + T^{23} + T^{24} \right), 
	\end{align} 
	where 
	\begin{align*}
	& T^{21} = e^{\lambda'}\int_{\mathcal{H}_u}\frac{1}{2}\sum_{j=1}^K (u_j - u_{0j})^2\alpha_{j}^{-1}
	e^{-\Phi_u^r(u)}\mu_r^u(du), \\
	& T^{22} = \frac{K}{2}|\lambda'|, \\
	& T^{23} = \int_{\mathcal{H}_{u}}\|Hu-d\|^2 e^{-\Phi_u^r(u)}\mu_r^u(du)
	\int_{\mathbb{R}}\frac{1}{2}e^{\tau'}e^{-\Phi_{\tau'}^r(\tau')}\mu_r^{\tau'}(d\tau'),  \\
	& T^{24} = \frac{N_d}{2}\int_{R}|\tau'|e^{-\Phi_{\tau'}^{r}(\tau')}\mu_{r}^{\tau'}(d\tau'). 
	\end{align*}
	Remembering that the operator $H$ is bounded and the specific forms of $a_{\tau'}(\epsilon,\tau')$ and $a_u(\epsilon, u)$, 
	we can obtain that the above four terms are all bounded. 
	The following inequality 
	\begin{align}\label{cond3}
	T^{3} := \sup_{\tau'\in T_{N}^{\tau'}}\sup_{\substack{\nu^{u}\in\mathcal{A}_{u}\\ \nu^{\lambda'}\in\mathcal{A}_{\lambda'}}}
	\int_{\mathcal{H}_u}\int_{\mathbb{R}}(\Phi^0 + \Phi)1_{A}(u,\lambda',\tau')\nu^{\lambda'}(d\lambda')\nu^{u}(du) < \infty
	\end{align}
	can be proved similarly, we omit the details.
	With the above calculations, we verified conditions (19) with $i=1,2,3$. 
	Now, we turn to verify conditions (20). 
	For conditions (20) with $i=2,3$, the inequalities could be verified similarly as for the case of $i=1$. 
	Hence, we only provide details when $i=1$ that is to prove
	\begin{align*}
	T^{4} := \sup_{\substack{\nu^{\lambda'}\in\mathcal{A}_{\lambda'}\\ \nu^{\tau'}\in\mathcal{A}_{\tau'}}} \!\!\int_{\mathcal{H}_u}\!\!\!\exp\bigg(\!\!-\!\int_{\mathbb{R}^2}(\Phi^0 + \Phi)1_{A^c}\nu^{\lambda'}\!(d\lambda')
	\nu^{\tau'}\!(d\tau')\!\bigg)\!\max(1,\|u\|_{\mathcal{H}_u}^2\!)\mu_{r}^u(du) < \infty. 
	\end{align*}
	Through a direct calculation, we find that 
	\begin{align*}
	-\int_{\mathbb{R}^2}(\Phi^0 + \Phi)1_{A^c}\nu^{\lambda'}\!(d\lambda') \nu^{\tau'}(d\tau') \leq &
	\frac{1}{2}\sum_{j=1}^K\alpha_{j}^{-1}(u_j-u_{0j})^2\int_{\mathbb{R}^{-}}(1-e^{\lambda'})\nu^{\lambda'}(d\lambda') \\
	& + \frac{K}{2}\int_{\mathbb{R}}|\lambda'|e^{-\Phi_{\lambda'}^{r}(\lambda')}\mu_r^{\lambda'}(d\lambda')  \\
	& + \frac{N_d}{2}\int_{\mathbb{R}}|\tau'|e^{-\Phi_{\tau'}^r(\tau')}\mu_{r}^{\tau'}(d\tau').
	\end{align*}
	Then we have 
	\begin{align*}
	T^{4} \leq C \int_{\mathcal{H}_u} \exp\bigg(
	\frac{1}{2}\sum_{j=1}^K\alpha_{j}^{-1}(u_j-u_{0j})^2\int_{\mathbb{R}^{-}}(1-e^{\lambda'})\nu^{\lambda'}(d\lambda') 
	\bigg) \max(1,\|u\|_{\mathcal{H}_u}^2)\mu_{r}^u(du).
	\end{align*}
	Considering $\int_{\mathbb{R}^{-}}(1-e^{\lambda'})\nu^{\lambda'}(d\lambda') < 1$ and the definition of $\mu_{r}^u$, 
	we know that the right hand side of the above inequality is bounded which completes the proof. 
	
	\section*{Acknowledgments}
	The authors would like to thank the anonymous
	referees for their comments and suggestions, which helped to improve the paper significantly. We also thank Ms. Ying Feng for her thorough polishing of this paper.
	This work was partially supported by the NSFC under the grants No. 11871392,
	and the National Science and Technology Major Project under grant Nos. 2016ZX05024-001-007 and 2017ZX05069.
	
	\bibliographystyle{plain}
	\bibliography{references}

\end{document}